\documentclass[12pt]{amsart}

\usepackage[dvipsnames]{xcolor}
\usepackage{latexsym}
\usepackage{tikz, pgfplots}

\usepackage{tikz}
\usepackage{tikz}
\usepackage{tikz-3dplot}
\usetikzlibrary{backgrounds}
\usepackage[dvipsnames]{xcolor}
\usepackage{subcaption}
\usepackage{animate}
\usepackage{pgfplots}\pgfplotsset{compat=1.13}
\usetikzlibrary{arrows, positioning, calc, patterns}
\usepackage{graphicx}

\usepackage[all]{xypic}
\usepackage{hyperref}

\usepackage{stmaryrd}

\usepackage{tkz-fct}  
\definecolor{bl}{RGB}{47,49,49}
\definecolor{skBlue}{RGB}{138 186 211}
\definecolor{egg}{RGB}{248,241,229}
\definecolor{citr}{RGB}{249,186,50}

\usepackage[papersize={210mm, 297mm}]{geometry}

\renewenvironment{frame}{}{}

\newcommand{\befr}{\begin{frame}}
\newcommand{\efr}{\end{frame}}
\usepackage{amsfonts,amssymb} 
\usepackage{stmaryrd}
\usepackage{amsmath}
\usepackage{mathtools}

\usepackage{cancel} 

\newtheorem{theorem}{Theorem}[section]
\newtheorem{lemma}[theorem]{Lemma}
\newtheorem{proposition}[theorem]{Proposition}
\newtheorem{remark}[theorem]{Remark}
\theoremstyle{definition}
\newtheorem{definition}[theorem]{Definition}
\newtheorem{example}[theorem]{Example}

\newtheorem{cor}[theorem]{Corollary}
\theoremstyle{remark}

\usepackage{xcolor}
\usepackage{xparse}
\usepackage{tikz}
\usetikzlibrary{matrix}
\usetikzlibrary{cd}

\newcommand{\ba}{\mathbf{a}}

\definecolor{amber}{rgb}{1.0, 0.49, 0.0}
\definecolor{amethyst}{rgb}{0.45, 0.31, 0.59}
\definecolor{applegreen}{rgb}{0.55, 0.71, 0.0}
\definecolor{asparagus}{rgb}{0.0, 0.42, 0.24}
\definecolor{darkbyzantium}{rgb}{0.36, 0.22, 0.33}
\definecolor{darklavender}{rgb}{0.45, 0.31, 0.59}
\definecolor{pr}{rgb}{1.0000000000000000000000000000000000000, 0.41176470588235294117647058823529411765, 0.38039215686274509803921568627450980392}

\definecolor{vp}{rgb}{0.050980392156862745098039215686274509804, 0.59607843137254901960784313725490196079, 0.72941176470588235294117647058823529412}

\definecolor{or}{rgb}{1.0000000000000000000000000000000000000, 0.64705882352941176470588235294117647059, 0}

\def\appendix#1{
\addtocounter{section}{1} \setcounter{equation}{0}
\renewcommand{\thesection}{\Alph{section}}
\section*{Appendix \thesection\protect\indent\quad
#1}
}
\begin{document}

\title[ $2-$valued groups, and Kontsevich-type polynomials, I]
{Algebraic $2-$valued group structures on $\mathbb P^1$, Kontsevich-type polynomials, and multiplication formulas, I }

\author{Victor Buchstaber, Ilia Gaiur and Vladimir Rubtsov}

\maketitle

\begin{flushright}{\em To 75-th anniversary of Andrei Bolibrukh}
\end{flushright}
\begin{abstract}
  The theory of a two-valued algebraic group structure on a complex plane and complex projective line is developed. In this theory, depending on the choice of the neutral element, the local multiplication law is given by the Buchstaber polynomial or the generalized Kontsevich polynomial. One of the most exciting results of our studies is a simple construction of a two-valued algebraic group on $\mathbb C$ different from known coset-construction.
\end{abstract}

        

\section{Dedication}
{\it A brilliant mathematician and an outstanding person, Andrei Bolibrukh was always held in great affection and high esteem by each of the three authors. It is well known that A. Bolibrukh's development as a mathematician began with algebraic topology, within the framework of Moscow seminars, where he first met and formed a friendship with the two senior authors. 

V.B. vividly recalls his presentations at M. Postnikov's seminar and the subsequent discussions with Andrei. Key aspects of the theory of 2-valued groups, to which certain parts of this work are devoted, turned out to be connected with branches of mathematics that had always been at the core of Andrei's interests. Later, Andrei became a friend and colleague at Steklov Mathematical Institute.

At the same time, for V.R., Andrei became (and remained) a senior colleague who was always ready to help, share knowledge, and explain challenging topics to a beginner. To a great extent, it was thanks to Andrei that he stayed in mathematics at the very start of his journey. This friendship endured throughout their lives in Moscow, Strasbourg, and Angers.

For the second author (I.G.), A.A. Bolibrukh became an iconic figure, a true Teacher and Mentor whose works on the analytic theory of differential equations, the Riemann-Hilbert problem, and the monodromy of solutions defined and shaped his primary mathematical interests, as well as the themes of his early works and Ph.D. thesis.

It is a great honor to dedicate this work to the memorable anniversary of Andrei Andreevich Bolibrukh.}

\section{Introduction}

The general concepts associated with multivalued groups had appeared in many contexts and have been intensively studied. The first author and S.P. Novikov in 1971 \cite{BuchNov}, came up with a product construction, with the origin in the theory of characteristic classes. This work initiated and gave birth to the theory of $N-$valued groups which was lntroduced and fully developed later by the first author.  Even for $N=2$ and in the case of finite groups their multivalued counterparts have quite different properties \cite{behravesh2012note}.

At this point, theories of formal, finite, finitely generated, algebraic, topological, and other multi-valued groups have been developed. The full account and survey of the theory of multivalued groups
was done in the paper \cite{Buch06}. 
 
In the present paper, we extend a theory of algebraic and topological two-valued groups on the complex projective line and study its applications. It is shown that some polynomial families  – {\it Buchstaber polynomials},  {\it generalized Kontsevich polynomials} and (as a special case) Buchstaber-Rees polynomials give rise to two-valued multiplication laws with identity element given by various points of a complex projective line for different families. These two 3-parametric families of symmetric polynomials in three variables are the central object of our interest in this paper.

\smallskip

 \subsection{Buchstaber polynomials}The family $B_{a_1,a_2,a_3}(x,y,z)$ which we call {\it Buchstaber polynomials} had appeared in 1990 \cite{Buch90}. The first author has classified the two-valued algebraic groups coming from the square modulus construction for formal groups with the multiplication law suggested by the Baker–Akhiezer elliptic function addition theorem. It reads as discriminant type polynomial family depending in 3 parameters $(a_1 ,a_2, a_3):$
\begin{equation}\label{buchstpol}
B_{a_1,a_2,a_3}(x,y,z):= (\sigma_1 - a_2\sigma_3)^2 -4(1+ a_3 \sigma_3)(\sigma_2 + a_1 \sigma_3),
\end{equation}
where $\sigma_k, k=1,2,3$ are elementary symmetric functions in $x,y,z.$

An explicit form of the Buchstaber polynomial family reads as
\begin{equation}\label{buchstpol1}
B_{a_1,a_2,a_3}(x,y,z) = (x+y+z - a_2x y z)^2 -4(1+ a_3 x y z)(xy + yz+ zx + a_1 xyz).
\end{equation}
 Zero locus $B_{a_1,a_2,a_3}(x,y,z)=0$  of this family defines a family of sextic hypersurfaces in $\mathbb C^3$
 
One can pose a question: for which values of parameters $a_i, i=1,2,3$  the family \eqref{buchstpol1} encodes the 2-valued group law in the coordinate origin of $0$ in $\mathbb C$? The answer is  brief and quite simple: this is true for {\bf any} $a_i$. 

Indeed, 
\begin{itemize}
\item take $y=0.$ Then from \eqref{buchstpol} we immediately get that $(z-x)^2=0$ and $0$ is the "neutral" of the supposed 2-value group law.
\item We know from the 2-valued group axiomatics (which we remind below in ...)  that the associativity $x\ast(y\ast z) = (x\ast y)\ast z$ one can interpret as equality of two polynomials in $x,y,z$ with coefficients depending on $a_1,a_2, a_3.$
\item The natural question of an extension (a "compactification") of the 2-valued group law with parameters $a_i$  from $\mathbb C$ to $\mathbb P^1$ leads to an independent appearance of elliptic
curves/functions and also demands the constraint of regularity. More precisely, the first author discovered that the polynomial \eqref{buchstpol} encodes the 2-valued group structure concerning the triple of variables $(x,y,z)$  in case when parameters $(a_1,a_2,a_3)$ belongs to an open everywhere dense subset in $\mathbb C^3$ where these parameters can are expressed via the Weierstrass $\wp-$ function with modular parameters $g_2, g_3$ in the following way ("uniformization"). Let $\alpha$ be a point of the elliptic curve $\mathcal E: v^2 = 4u^3 - g_2u - g_3$, then polynomial \eqref{buchstpol} defines the 2-valued group  in the case when the following equalities hold
\begin{equation}\label{unipar}
a_1 = 3\wp (\alpha), \quad a_2 = 3\wp(\alpha)^2 - g_2/4,\quad  a_3 = 1/4(4\wp(\alpha)^3 - g_2\wp(\alpha) -g_3).
\end{equation}
Using this uniformization we show that the algebraic law on $\mathbb{C}$ is defined for any choice of $a_1, a_2$ and $a_3$. We also show that this algebraic law extends to $\mathbb{P}^1$ when $a_1,a_2$ and $a_3$ are in the complement to some discriminantal locus (see theorem \ref{basiclaw} for the details).

It gives a universal parameter uniformization in the class of algebraic $2$-valued group appeared from the "square modules construction" and may be viewed as a coset group $\mathcal E/\imath$, where $\mathcal E$ is the elliptic curve (considered as an abelian group) and $\imath: \mathcal E\to \mathcal E$ is the antipodal involution which sends $v$ to $-v$. 

Another interesting interpretation of the universal uniformization can be described in terms of {\it groupoid operation} which was a partial case of hyperelliptic curve jacobians addition laws which were studied in terms of Weierstrass $\sigma-$functions in \cite{buchstaber2005addition}. We discuss it briefly in subsection 4.3.
\end{itemize}

The above-mentioned algebraic 2-valued groups define a quadratic addition law.  We show that there is a reversible transformation of the corresponding Buchstaber polynomial into the so-called generalized Kontsevich polynomial - the second object of our interest and study.

\subsection{Kontsevich polynomials and explicit Geometric Langlands computations}The first appearance of these polynomials, up to our knowledge, is due to the explicit computations in {\it Geometric Langlands Correspondence.} The Langlands correspondence relates Galois representations to automorphic
data connecting Frobenius eigenvalues on the Galois side to the Hecke eigenvalues on the automorphic side.
In his famous paper \cite{Konts2007} M. Kontsevich has proposed a computation to make the Geometric Langlands Correspondence explicit for $SL_{2}-$ local systems on $\mathbb P^1\setminus\{0, t, 1, \infty\}$ with unipotent local monodromies. Automorphic forms in his description are complex-valued functions of the finite field  $\mathbb F_q,\, {\rm char}(\mathbb F_q) \neq 2.$

Let us describe some nodes of his construction in the context of our paper.
Kontsevich introduced an associative and commutative algebra structure on the vector space of $q\times q$ matrix $T_x$ (for each $x\in \mathbb F_q$)   whose entries at $(y,z)$ are functions  of root numbers (in $\mathbb F_q$) of the universal discriminant type polynomial with integral coefficients 
\begin{equation}\label{kontspol}
u^2 - P_t (x,y,z) = u^2 -[(xy+ yz+zx-t)^2 - 4xyz(1+ t+ (x+y+z))].
\end{equation}
More precisely, the  matrix coefficients are computed via the following formula:

If $x\neq \infty,y\neq \infty$ and the $z-$discriminant of $P_t (x,y,z)$ is non=zero, then
\begin{equation}\label{Hecke1}
(T_x)_{yz} = \begin{cases}
2-\sharp\{u\in \mathbb F_q | u^2 = P_t(x,y,z)\} & z\neq \infty \\
2-\delta_{x,y} & z=\infty .
\end{cases}
\end{equation}
In the remaining cases matrix element reads
\begin{equation}\label{Hecke2}
(T_x)_{yz} = \begin{cases}
1-q & \text{if} \, z =\frac{(xy+t)((x+y)+2(1+t)xy}{(x-y)^2} \\
2 & \text{otherwise}.
\end{cases}
\end{equation}
Matrices $T_x\in {\rm Mat}_{q\times q}(\mathbb F_q)$ may be considered as {\it Hecke operators}. These matrices commute for different points, i.e. $[T_x, T_y]=0$. Moreover, they form an associative algebra which reads as
\begin{equation}\label{assprod1}
T_x·T_y=\sum_{z\in F_q} C_{xyz}T_z, {\rm where}\quad  C_{xyz}: =(T_x)_{yz}.
\end{equation}

Similar considerations and computations were also suggested by T. Thomas in his Master Thesis Project \cite{Thomas}, where, in particular, he had checked carefully the associativity of the product \eqref{assprod1}. Also in the dissertation of  N. de Bos \cite{deBos} the same problem was considered and extended.

\subsection{Kontsevich polynomials and elliptic curve group law} Another appearance of the Kontsevich polynomial is due to the multiplication formula for the Heun equations \cite{MKSO}. However, it may be obtained as a elliptic curve group law projected on the base. In particular, the polynomial family in the RHS of \eqref{kontspol} (the "usual" Kontsevich polynomial family) is a $t-$discriminant $D_{\lambda}(x,y,z)$ of  $P_{t}(x,y,z).$ 
Moreover, the Kontsevich polynomial is symmetric in $x,y,z.$ Indeed, the following formula holds 
\[D_{\lambda}(x,y,z)  = \operatorname{disc}_t (f(t) - (t-x)(t-y)(t-z))
\]
where $f(t) = t(t-1)(t-\lambda) \in \mathbb C[t]$ and $\lambda \in  \mathbb{C}^1$.

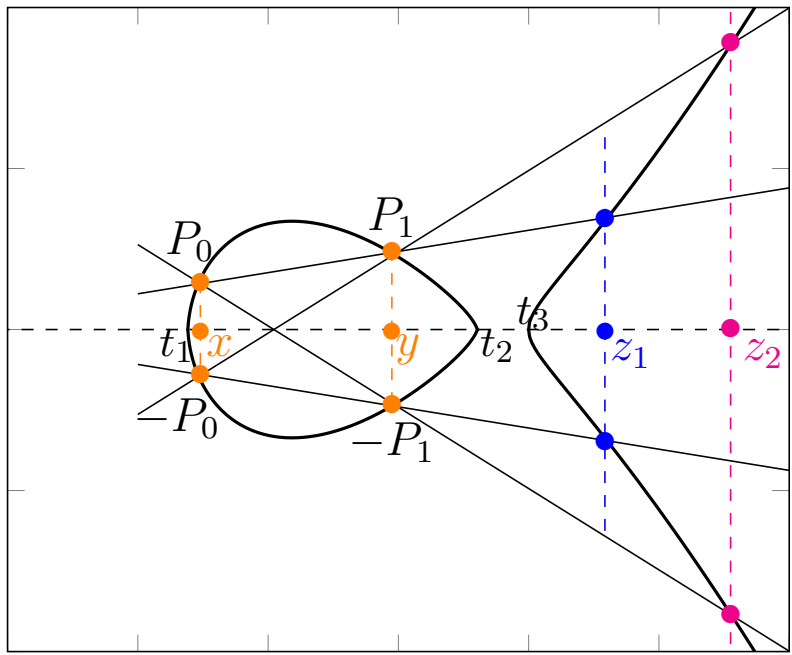
\begin{figure}
  \centering
  
         \begin{tikzpicture}[scale=1.5]
\begin{axis}[
            xmin=-3.0,
            xmax=3,
            ymin=-4,
            ymax=4,
            domain=-1.521379707:2,
            samples=500,
            smooth,
            yticklabels={,,},
            xticklabels={0},
        ]
       
            \addplot[orange, domain=1.0:-1.0, dashed]({-0.05},{x});
          \addplot[orange, domain=0.6:-0.6, dashed]({-1.52},{x});
         
          \addplot[dashed, domain=-4:5]{0};
        \addplot[domain=-1.7:0.61,thick]{sqrt(x^3-2*x+1)-0.1};
       
        \addplot[domain=-1.7:0.61,thick]{-sqrt(x^3-2*x+1)+0.1};
        \addplot[domain=-2:3]{1/3.8*x+0.97};
        \addplot[domain=-2:3]{-1/3.8*x-0.96};
        \addplot[domain=-2:3]{1.01*(x+1.52)-0.57};
        \addplot[domain=-2:3]{-1.01*(x+1.52)+0.57};
        \node at (-1.7,-0.2) {\small $t_1$};
        \node[anchor=west] at (0.5,-0.2) {\small $t_2$};
        \node[anchor=east] at (1.3,0.2) {\small $t_3$};
       
        \addplot[domain=1:4,thick]{sqrt(x^3-2*x+1)};
        \addplot[domain=1:4,thick]{-sqrt(x^3-2*x+1)};
        \addplot[blue,domain=-2.5:2.5,dashed]({1.585351681},{x});
        \addplot[magenta,domain=-4.5:4.5,dashed]({2.55},{x});
        \node [magenta] at (2.55,3.55)
        {\textbullet};
        \node [magenta] at (2.55,-3.55)
        {\textbullet};
        \node [magenta] at (2.55,0)
        {\textbullet};
        \node [magenta] at (2.8,-0.3)
        {$z_2$};
       
        \node [blue] at (1.585351681,-1.4)
        {\textbullet};
        \node [blue] at (1.585351681,1.37)
        {\textbullet};
        \node [orange] at (-0.05,0.95)
        {\textbullet};
        \node [orange] at (-0.05,-0.95)
        {\textbullet};
        \node  at (-0.05,1.4)
        {$P_1$};
        \node  at (-0.05,-1.4)
        {$-P_1$};
        \node  at (-1.6,1.1)
        {$P_0$};
        \node  at (-1.7,-1.1)
        {$-P_0$};

        \node [orange] at (-0.05,-0.035) {\small \textbullet};
        \node [blue] at (1.585351681,-0.035) {\small \textbullet};
        \node [orange] at (-1.52,-0.035) {\small \textbullet};
        \node[anchor=west,orange]  at (-0.15,-0.25)
        {$y$};
        \node[anchor=west,orange]  at (-1.6,-0.2)
        {$x$};
        \node[anchor=west,blue]  at (1.5,-.3)
        {$z_1$};
       
        \node [orange] at (-1.52,0.57) {\textbullet};
        \node [orange] at (-1.52,-0.57)
        {\textbullet};
        \end{axis}
\end{tikzpicture}

\caption{Geometric construction for the addition law on elliptic curve $y^2 = (t - t_1)(t - t_2)(t - t_3)$
with a choice of neutral element at $\infty$ and corresponded coset 2-valued group law $x*y = [z_1,z_2]$ on $\mathbb P^1.$}
\label{fig:elllaw}
\end{figure}

In this paper, following to \cite{GMRvS} we study a similar discriminant family for a general cubic polynomial
$$f(t) = t^3 + at^2 + bt + c
$$which we call {\em generalized Kontsevich polynomial} which is given by
\begin{equation}\label{eq:disc1}
    D_{a,b,c}(x,y,z) = \operatorname{disc}_t(t^3 + at^2 + bt + c -  (t-x)(t-y)(t-z)).
\end{equation}
and more precisely given by
   \begin{align}\label{genkontspol}
     D_{a,b,c}(x,y,z)=\left(x -y \right)^{2} z^{2}-2 \left(2 a x y + \left(x +y \right)( x y + b) +2 c \right)z+\nonumber \\+(yx-b)^2 -4 c \left(x +y+a\right).
\end{align}
The corresponding loci represent an image of an Abelian law for the corresponding elliptic curve 
$$
\mathcal{E}: \quad w^2 = t^3 + at^2 + bt + c,
$$
projected to the base curve. The construction goes as follows:  let $x$ and $y$ be the points in $\mathbb{P}^1$.  Consider elliptic curve $\mathcal{E}$ as a double cover of $\mathbb{P}^1$ ramified at the roots of $f(t)$ and $\infty$. Let $x_{+}, x_{-},y_{+}$ and $y_{-}$ be the points on $\mathcal{E}$ which are pre-images of the projection map $\pi: \, \mathcal{E}\,\rightarrow \,\mathbb{P}^1$ for $x$ and $y$ correspondingly. Take all possible sums $x_{\pm}\,+\,y_{\pm} \in \mathcal{E}$  on the elliptic curve $\mathcal{E}$. Roots $z_{1,2}$  of \eqref{genkontspol} are projections of $x_{\pm}\,"+"\,y_{\pm}$ to the base line $\mathbb P^1$ (dashed line on the picture). This construction is illustrated in Figure \ref{fig:elllaw}. Described construction is an example of more general coset construction from multi-valued groups theory.

%


\begin{remark}
 The generalized Kontsevich polynomial \eqref{genkontspol} is $\mathfrak S_3-$ invariant under the cyclic permutations of $x,y$ and $z-$ variables.
 One can write it directly in terms of elementary symmetric functions as:
 \begin{equation}\label{genkontsimm}
D_{a,b,c}(\sigma_1,\sigma_2,\sigma_3)=(\sigma_2 - b)^2 - 4 (\sigma_3 + c) (\sigma_1+a).
 \end{equation}
The symmetric form of Kontsevich polynomials will be used in subsection 4.2.
\end{remark}


\subsection{Structure of the paper}  The paper contains the Introduction and three sections. The detailed Introduction advertises all the main "actors" of the article. Then in Section \ref{section:2gr} we propose a few comprehensive reminders of the 2-valued group theory covering both formal and algebraic constructions. Some important examples are discussed in detail.
 
The first part of Section \ref{section:P1} contains a description of the Buchstaber polynomial \eqref{buchstpol1} appearance in the context of the "square modulus construction". We demonstrate of the coefficient rationality for defining quadratic equation for 2-valued group operation in full detail (Th. \ref{th:rat})  using both direct arguments and an elegant approach proposed in the Burnside classical notes.

The Burnside determinant triviality (which is a "0-case" of the classical Frobenius-Stickelberger formula) implies the existence of  a link between the special case (when the parameter $a_1=0$) of \eqref{buchstpol1} and the Kontsevich polynomial which reads off the addition law on the "bare" elliptic curve.
 
The second part of the section occupies two important and interesting questions: how the generalized Kontsevich polynomial naturally appears as a result of the Weierstrass addition law deformation and application of a non-triviality of Frobenius-Stickelberger formula. Another question (subsect. \ref{subs:alg}, th. \ref{basiclaw}) was mentioned in the introduction: to give a new description of "extendable" 2-group laws on $\mathbb P^1$. The classification of such laws was known due to Buchstaber and Rees. Still, we are clarifying the appearance and character of obstructions to an extension of  (non-obstructed) 3-parametric 2-valued group laws on $\mathbb C$ to a 2-valued group structure on $\mathbb P^1$. The Th. \ref{basiclaw} fully describes 2-valued "extendable" from $\mathbb C$ to $\mathbb P^1$ 2-valued groups defined by the Buchstaber polynomial. 

We provide finally in subsection \ref{subs:moeb} a few examples of "deformed" $\mathbb P^1-$2-valued group structure given by the M\"{o}bius action on the elementary quadratic 2- law equation. We remark that the first study of the M\"{o}bius action on quadratic symmetric polynomials with three variables was studied by V.Dragovi\v{c} \cite{Drag}.

As we have already seen in the examples of Section \ref{section:P1}, there exists a birational automorphism $\mathbb P^2$ transforming the Kontsevich polynomial into the Buchstaber polynomial (up to the Jacobian factor/divisor. In Section \ref{sec:BK} we describe this correspondence explicitly drawing attention to the geometric picture of the story. Our main (rather straightforwardly computational) result is summarized in the Th. \ref{th:KB1}. Then we briefly remark the {\it discriminantly-splits} of both polynomial families and prove that {\it all} symmetric discriminantly–split polynomials with three variables, such that they have degree two in each of them and satisfy criteria of 2-valued group law coincide with Buchstaber-Kontsevich type polynomials. We describe also an interesting algebraic and geometric property of the intersection loci for two families of polynomials and establish a link with Beauville's classification of elliptic modular surface fibrations.

The subsection \ref{subs:var} stands a little bit apart but we have included it because of transparent and computationally simple arguments of the correspondence between Buchstaber and Kontsevich polynomials after the passage to {\it $\star-$involution} map of (localized) in $\sigma_3$ Laurent polynomials in thre elementary $\sigma_i, i=1,2,3$ coordinates. The theorem \ref{th:KB1} re-written after this involution becomes practically tautologic (see Th. \ref{th:var}). 

We discuss briefly the content of the second part of our paper in Section \ref{sec:concl}.

 {\subsection*{Acknowledgements} V. R. and I.G. greatly acknowledged the IHES hospitality during two visits in 2023 and 2024, where this paper was started and finished. We are much indebted to D. Gurevich and D. Talalaev for numerous discussions of multivalued groups and $n-$coalgebra notions. V.R. work was in part supported by the Project RSF Grant No. 23-41-00049.  He also acknowledges the France 2030 framework program Centre Henri Lebesgue ANR-11-LABX-0020-01 for creating an attractive mathematical environment.}

\section{Two–valued groups in all its states}\label{section:2gr}

\subsection{Two-valued group generalities.}

We remind for completeness of the basic definitions of two-valued group theory (based for example on \cite{Buch06, BGV})

A {\it two-valued group} is a set $X$ with a two-valued multiplication law $\ast : X\times X \to {\rm Sym}^2 (X),$ {\it identity element} $e\in X,$ and the {\it inverse. map} 
$x\to x^{-1}$ with the following properties:
\begin{itemize}
\item {\it Associativity:} for any three elements $x,y,z \in X$ there is a coincidence of four- element multisets

$$ 
(x\ast y)\ast z = x\ast (y\ast z)
$$
\item {\it Strong identity:} for any element $x$ there are equalities

$$
e \ast x = x \ast e =  [x, x].
$$
\item {\it Existence of the inverse:} for any element $x \in X,$ each of the multisets $x \ast x^{-1}$ and $x^{-1}\ast x$ contains the identity $e.$ 
\end{itemize}
\begin{remark}
The above axioms demand an inverse element's existence only and do not imply its uniqueness. Thus,  the uniqueness of the inverse is an additional strong condition on a 2-valued group structure. For example, the coset construction  (see below) satisfies this condition.
\end{remark}

A two-valued group is called {\it commutative} if $x\ast y=y\ast x$ for all $x,y \in X.$ We shall call an element $x$ of a two-valued group $X$ {\it a weak involutive} if $x^{-1} = x$ or, equivalently, if the multiset $x\ast x$ contains the identity $e.$ 
An element $x\in X$ is a {\it strong involutive} if $x\ast x = [e, e].$ A two-valued group $X$ is {\it involutive} if it consists of weak involutions, that is, 
$x^{-1} = x$ for all $x \in X.$
Classification of 2-valued involutive commutative groups was done in \cite{BGV}. Later in \cite{Gaif} it was proven that every 2-valued involutive group is commutative, which leads to a full classification of 2-valued involutive groups.

\subsection{Algebraic families of 2-valued group laws on $\mathbb C$}
What we shall understand under an {\it algebraic family of 2-valued group laws on $\mathbb C$?}

\begin{definition}
Let $\bar a = (a_1,\ldots, a_n)$ be a {\it vector of complex parameters}. Consider three polynomial families $A_{\bar a}(x,y), B_{\bar a}(x,y), C_{\bar a}(x,y)\in \mathbb C[x,y]$.
Then the family 
\[\mu_{\bar a} : \mathbb C\times \mathbb C \to {\rm Sym^2 \mathbb C}\]
of algebraic mappings defined by $\mu_{\bar a} (x,y) = [z_1,z_2],$ where $z_{1,2}$ are roots of $A_{\bar a}z^2 +B_{\bar a}z +  C_{\bar a} =0,$
gives an algebraic family of 2-valued group laws on $\mathbb C$ if
\begin{itemize}
\item a (strong) identity element $e$ is defined.  
\item for any $x\in \mathbb C$ exist $ \nu(x)\in \mathbb C$ such that  $\mu_{\bar a}(x,\nu(x))$ contains $e$.
\item {\it Associativity:} for three complex numbers $x,y,z \in \mathbb C$ there exist another triple $(u,v,w) \in \mathbb C$ such that the following associativity condition takes place:

let a pair $(X_1,X_2)$ be the pair of roots for the quadratic equation
 \[ A_{\bar a}(x,y)X^2 +B_{\bar a}(x,y)X +  C_{\bar a}(x,y)= 0\]
 while a pair $(Y_1,Y_2)$ be similarly the pair of roots for 
 \[A_{\bar a}(y,z)Y^2 +B_{\bar a}(y,z)Y +  C_{\bar a}(y,z)= 0.\]
 Then the system of quadratic equations:
\begin{equation*}\label{assfamcorrespL}
\begin{cases}
\centering
 &A_{\bar a}(X_1, z)u^2 +B_{\bar a}(X_1,z)u +  C_{\bar a}(X_1,z)= 0\\
 &A_{\bar a}(X_2, z)u^2 +B_{\bar a}(X_2,z)u +  C_{\bar a}(X_2,z)= 0 
\end{cases}
\end{equation*}
considering as the $4-$th order polynomial equation

\begin{align}\label{eq:11}
\begin{split}
A_{\bar a}(X_1,z)A_{\bar a}(X_2,z)u^4 +\left[A_{\bar a}(X_1,z)B_{\bar a}(X_2,z) +B_{\bar a}(X_1,z)A_{\bar a}(X_2,z)\right]u^3 +\\
+ \left[ A_{\bar a}(X_1,z)C_{\bar a}(X_2,z)+C_{\bar a}(X_1,z)A_{\bar a}(X_2,z) +B_{\bar a}(X_1,z)B_{\bar a}(X_2,z)\right]u^2 +\\
+ \left[B_{\bar a}(X_1,z)C_{\bar a}(X_2,z)+C_{\bar a}(X_1,z)B_{\bar a}(X_2,z)\right]u + C_{\bar a}(X_1,z)C_{\bar a}(X_2,z) =0.
\end{split}
\end{align}

gives the same  4-set of roots $(x, y, z,w)\in \mathbb C^4$ as the system
\begin{equation*}\label{assfamcorrespR}
\begin{cases}
\centering
&A_{\bar a}(x,Y_1)v^2 +B_{\bar a}(x,Y_1)v +  C_{\bar a}(x,Y_1)= 0\\
&A_{\bar a}(x, Y_2)v^2 +B_{\bar a}(x,Y_2)v +  C_{\bar a}(x,Y_2)= 0.
\end{cases}
\end{equation*}
considering as the $4-$order polynomial equation

\begin{align}\label{eq:12}
\begin{split}
A_{\bar a}(x,Y_1)A_{\bar a}(x,Y_2)v^4 +\left[A_{\bar a}(x,Y_1)B_{\bar a}(x,Y_2) +B_{\bar a}(x,Y_1)A_{\bar a}(x, Y_2)\right]v^3+\\
+ \left[ A_{\bar a}(x,Y_1)C_{\bar a}(x,Y_2)+C_{\bar a}(x,Y_1)A_{\bar a}(x, Y_2) +B_{\bar a}(x,Y_1)B_{\bar a}(x,Y_2) \right]v^2 +\\
+ \left[B_{\bar a}(x,Y_1)C_{\bar a}(x,Y_2)+C_{\bar a}(x,Y_1)B_{\bar a}(x, Y_2)\right]v + C_{\bar a}(x,Y_1)C_{\bar a}(x,Y_2) =0
\end{split}
\end{align}
\end{itemize}
\end{definition}

\begin{proposition}\label{prop:3.3}
A verification of the last axiom (associativity) is reduced to an algebraic problem.
\end{proposition}
\begin{proof} 
We observe that both systems \eqref{assfamcorrespL} and s \eqref{assfamcorrespR} are symmetric with respect to the "flips" $X_1 \leftrightarrow X_2$ and 
$Y_1 \leftrightarrow Y_2$ respectively. Hence, there exists rational  $\Theta_1(x,y) =\frac{B_{\bar a}}{A_{\bar a}}$ and  $\Theta_2(x,y)=\frac{C_{\bar a}}{A_{\bar a}}$  and 
the associativity conditions I-IV from ch. 2 in \cite{Buch2} (which were formulated there for formal series $\mathbb C[[x,y,z]]$) are equivalent to the hypothesis (3). 
But the multiplication by common denominators  $A(X_1,z)A(X_2,z)$ and  $A(x, Y_1)A(x, Y_2)$ makes the conditions (12-14) and (15-17) polynomial.
\end{proof}

We shall describe in section 4 a geometric interpretation of these conditions on the mapping family $\mu_{\bar a}$ when it gives an algebraic 2-valued group law
and prove that the corresponding polynomial family belongs to Buchstaber–Kontsevich generalized polynomials.

\subsection{2-coset construction}



  Let $G$ be any group and ${\rm A}(G)$  a subgroup of ${\rm Aut}(G)$, such that $\natural(A)=2$. It means that $A$ is given by involution and identity map.
  Then one can define a two-valued group structure on $X = G/A$ as follows: let $\pi  : G \to G/A$ be the quotient map and define $\ast : X\times X \to {\rm Sym}^2 (X)$ by

$$ x\ast y = \pi(\pi^{-1} (x)\circ \pi^{-1} (y))
$$ 
where $\circ$ denotes the multiplication on $G$, extended to sets.
A 2-valued group of this type will be called a {\it 2-coset group.}

If $A = \langle \operatorname{id},a\rangle$, Then the elements of $X$ can be written as
the pairs $(g, a\dot g)$and 
$$(g, a(g)) \ast  (h, a(h))= [(gh, a(g)a(h)), (ga(h),a(g)h)].
$$One may easily check that this construction satisfies all the axioms of $2$-valued groups, 
with inverse given by $\operatorname{inv}(x)=\pi(g^{-1}),$ where $\pi(g) =x.$

%

\subsection{2-valued group. "Toy model" polynomial example}

 Consider two important examples of 2-valued groups. The set $X=\mathbb C$ and we shall identify $\mathbb C^2$ with ${\rm Sym}^2 (\mathbb C)) = \mathbb C^2/{\mathfrak{S_2}}= \mathbb C^2/{\mathbb Z_2}$ via 
 
$$\mathcal S: \mathbb C^2 \to \mathbb C^2,\quad \mathcal S(z_1,z_2) =  (e_1(z_1, z_2), e_2(z_1, z_2)),
$$
where $e_k$ are elementary symmetric polynomials of degree $k=1,2$ in $z_{1,2}$. It is often convenient to write the map $\mathcal S$ as the quadratic polynomial $z^2 - e_1 z + e_2$  whose roots are $[z_1, z_2].$ 

\begin{example}\label{basic-2}
Consider $(\mathbb C, +)$ with its automorphism $z \to - z.$
We obtain an $2-$valued group structure on $\mathbb C/\mathbb Z_2.$ But $\mathbb C/\mathbb Z_2$ is homeomorphic to $\mathbb C$ by the mapping 
$ \mathbb C \to \mathbb C/\mathbb Z_2$(2-branched  covering ramified at $z=0$) defined by $\pm z \to z^2.$
The multiplication is then given, for $x,y \in \mathbb C$ by

$$ 
x\ast y = [(\sqrt x +\sqrt y)^2, (\sqrt x -\sqrt y)^2].
$$
The unit $e$ is $0$ and the inverse of $z$ is $z.$  Using the map $\mathcal S$, one obtains the following polynomial law $p_2(x,y,z):$

\begin{align}\label{twoval}
\begin{split}
p_2(x,y,z) = z^2 - 2(x+y) z  + (x-y)^2 &= \\
(z +x + y)^2 - 4(xy + yz + zx) &=\\
{\rm disc}_t [(t^2 + (z + x + y)t + (xy + yz + yz)].
\end{split}
\end{align}

The $2-$valued multiplication is described by this polynomial (which is a special partial case of the Buchstaber polynomial and corresponds to zero parameter choice) $B_{0,0,0}(x,y,z) = p_2(z; x, y) .$

$$
p_2(z; x, y) = (z- ((x)^{-1}\ast (y)^{-1})_1)(z - ((x)^{-1}\ast (y)^{-1})_2) =(z -(\sqrt x-\sqrt y)^2)(z -(\sqrt x+\sqrt y)^2),
$$
whence the product $(x)^{-1}\ast (y)^{-1}= x\ast y$ is defined by $z-$roots of the equation $p_2 = 0.$
\end{example}

More general, similar result holds for the higher Buchstaber-Rees polynomials. Namely, one may consider a coset construction given by $\mathbb{C}/\mathbb{Z}_N$ via multiplication by the $N$-th root of unity. In that case $N$-valued group law reads
$$
x\star_N y = \left[\left(x^{1/N}+\chi^k y^{1/N}\right)^N,\quad k=0\dots N-1\right],\quad \chi^N = 1.
$$
It was proven in \cite{GRvS} that multiplication kernels for the higher Bessel equations have singularities exactly along Buchstaber-Rees polynomials which reads as
$$
p_N(z; x, y) = (z- ((x)^{-1}\star_N (y)^{-1})_1)(z - ((x)^{-1}\star_N (y)^{-1})_2) \dots (z - ((x)^{-1}\star_N (y)^{-1})_N).
$$
For more details of this case we refer to \cite{GRvS}.

\begin{example}\label{basic-3}
Another important multiplicative 2-valued group on $\mathbb C$ can be described as follows. Let $G = \mathbb C^{*} = (\mathbb C\setminus 0).$ The automorphism
$z \to z^{-1}$  generates the order 2 group automorphism $A$. The space $X = \mathbb C^{*}/A$ is identified with $\mathbb C$ by means of the map 
$\mathbb C^{*} \to \mathbb C$ given by $\pi(z) = \frac{1}{2} (z+z^{-1})$ and $\pi^{-1}(x) = x\pm \sqrt{(x^2 -1)}$
and we obtain
\begin{equation}
x\ast y =\left [xy + \sqrt{(x^2 -1)(y^2 -1)},\,  xy - \sqrt{(x^2 -1)(y^2 -1)} \right]:
\end{equation}
The right-hand side can be rewritten as the quadratic polynomial in $z:$
\begin{equation}
z^2 - 2xyz + x^2 + y^2  -1:
\end{equation}
The unit is 1, and each element's inverse is the element itself. Under the change of
variables 
$x \to x + 1; y\to y + 1; z\to z + 1,$ the polynomial becomes
\begin{equation}\label{mord}
z^2 - 2(x+y +xy)z + (x-y)^2.
\end{equation}
The latter polynomial is nothing but the Mordell modification of a generalized Markov cubic polynomial \cite{BuchVes}.
\end{example}

\subsection{Algebraic 2-valued groups}
This example is the first important example of involutive  {\it algebraic two-valued multiplication law}. Such a law  $(x, y) \mapsto z = x\ast y$ is given in local coordinates by the equation $F (x, y, z) = 0$, where $F (x, y, z)$ is a polynomial, which is symmetric in the variables $x, y, z$ and of degree 2 in each of them and we are fully in this framework. The two-valued multiplication law must satisfy the associativity condition

$$
(x\ast y)\ast z = x\ast (y\ast z)
$$
in the following sense:

 Let $(x, y, z)$ and $(u,v,w)$ be two sets of variables in $\mathbb C^3.$  The algebraic associativity condition means that two systems of equations

\begin{equation}\label{asscorresp}
\centering
\left\{\begin{split}
F(x, y, u) &= 0\\
F(u, z, w)&= 0 \\ 
\end{split}\right. \quad
{\rm and} \quad
\centering
\left\{\begin{split}
F(u, z, v) &= 0\\
F(x, v, w) &= 0 \\ 
\end{split}\right.
\end{equation}
forms a {\it correspondence} in $\mathbb{C}^6$: which means that they define (after exclusion of $u$ and $v$) the same quadruple of points $(x,y,z,w)\in \mathbb C^4.$
We remark also that $F(0,y,z) = (y- z)^2$ (and hence also $F(x,y,0) = (y- x)^2$ and $F(x,0,z) = (x- z)^2$). Such algebraic symmetric 2-valued formal groups are called {\it involutive.}  It is well–known that all such groups are given by the Buchstaber polynomials up to a 2-valued group  unity choice (\cite{BuchVes} and the discussion in our paper below.)

\subsection{Two-valued formal groups.}

We are briefly reviewing main definitions and notions specifying the {\it formal 2-valued groups} within the general setting of 2-valued group laws. We shall assume in this subsection that our basic commutative unital ring $R$ is a $\mathbb Q-$algebra. We shall consider {\it 2-valued formal series} $\mathcal F(x,y)$ over $R$ the following quadratic equation
\begin{equation}\label{twoser}
\mathcal F^2 - \Theta_1(x,y)\mathcal F +\Theta_2 (x,y) = 0, \quad \Theta_{1,2}\in R[[x,y]]
\end{equation}
We shall omit all details and precisions about solutions of such equations sending a reader to the original paper of the first author  \cite{Buch2}, and concentrate our attention on the following specialization of a two-valued {\it formal} group definition:
\begin{definition}\label{univ}
A one-dimensional commutative 2-valued formal group over $R$ is a 2-valued formal series $\mathcal F(x,y)\in R[[x,y]]$ given by the equation \eqref{twoser} such that 
\begin{itemize}
\item $\mathcal F(x,0)$ is a two-valued series 
$$\mathcal F^2 - 2x\mathcal F + x^2 = 0;
$$
\item The associativity correspondence (see \eqref{asscorresp}) holds
\item $\mathcal F(x,y) = \mathcal F(y,x).$
\end{itemize}
\end{definition}
We shall assume additionally that $R$ has no zero divisors. Let $X_1,X_2$ be the roots of \eqref{twoser}. Then $X_1 = f(x,y) \in R[x,y] \iff X_2=X_1$ 
and $f(x,y)$ is a formal commutative group over $R$. The equation \eqref{twoser} transforms to the equation
$$
\mathcal F^2 - 2f(x,y)\mathcal F + f^2 (x,y) = 0.
$$
Another important examples is given by the so-called elementary two-valued formal groups given by the equations
$$
\mathcal F^2 - 2(x+y)\mathcal F +  (x+y)^2 = 0,
$$
and
$$
\mathcal F^2 - 2(x+y)\mathcal F +  (x-y)^2 = 0.
$$
These groups are also called {\it two-valued formal groups of basic type.}

\subsection{Two-valued formal groups generalities}
Let suppose that the ring $R[[x]]$ is enabled with a {\it commutative coalgebra} structure with the comultiplication

$$\Delta_{x}^{y} : R[[x]]\mapsto R[[x]]{\hat \otimes}_{R}R[[x]] \simeq R[[x,y]],
$$
such that $0-$ two-side coalgebra co-unity ($\Delta_{x}^{y} (f)\vert_{y=0} =f(x), \Delta_{x}^{y} (f)\vert_{x=0} =f(y).$) This coalgebra will be a homogeneous if
the comultiplication maps degree $n$ monomials in degree $n$ homogeneous polynomials in both variables.

It was shown in \cite{Buch2} that for any two-valued series satisfying to \eqref{twoser} with $\Theta_{1,2}(0,0) =0$ the following comultiplication
$$
\Delta_{x}^{y} : R[\frac{1}{2}][[x]] \to R[\frac{1}{2}][[x,y]], \quad \Delta_{x}^{y} (x^k) = \frac{1}{2}(\mathcal F^{(1)})^k (x,y) +(\mathcal F^{(2)})^k (x,y)).
$$
where $\mathcal F^{(i)}, i=1,2$ are symbolic solution sets of the quadratic equation for the 2-valued formal series \eqref{twoser}.
To check the co-module structure given by the mapping it is enough to impose their linearity condition and the classical co-associativity axiom. 
\smallskip
One of important observations in \cite{Buch2} is the following Proposition (Lemmas 3.1 and 3.2 in \cite{Buch2}):
\begin{proposition}\label{propcoalg}
The quadratic equation \eqref{twoser} is a two-valued formal group if and only if  the homomorphism

$$\Delta_{x}^{y} : R[\frac{1}{2}][[x]] \to R[\frac{1}{2}][[x,y]]
$$
defines a structure of a commutative coalgebra in the ring $R[\frac{1}{2}][[x]].$ 
\end{proposition}

It is important to stress the difference in axiomatics of formal and algebraic 2-valued groups. In formal 2-valued groups, $0$ is a "unit" (neutral element) of the group and this determines the type of group law.
Whereas for algebraic 2-valued groups, we should postulate only the existence of a unit. In other words, we can claim that there is only a notion of an algebraic 2-valued group but there is not a definite law behind it. The law appears only if we choose the unity and only then the law can be written in a neighborhood of it.  
We shall illustrate this principle below.

\section{2-valued group structures on $\mathbb P^1.$}\label{section:P1}
\subsection{Buchstaber elliptic family of polynomials}
Buchstaber had classified the two-valued algebraic groups arising from the square modulus construction for formal groups with the multiplication law suggested by the addition theorem for Baker–Akhiezer elliptic functions \cite{Buch90}. Square modulus construction reads as follows. Consider polynomial
\[B(x,y,z):=z^2 - \Theta_1(x,y)z + \Theta_2(x,y)=0
\]
and the formal group $F(u,v)$ with exponent $$f(u) = \frac{\sigma(u)\sigma(\alpha)}{\sigma(\alpha -u)}\exp(-\zeta(\alpha)u), $$ where $\sigma(u)$ - Weierstrass sigma-function and $\zeta(u) = \frac{d\log(\sigma(u))}{du}.$

Set
\[
\Theta_1(x,y) = Z_{+} + Z_{-}; \quad   \Theta_2(x,y) = Z_{+}Z_{-},
\]
with
\[Z_{+} = F(u, v)F(\bar u, \bar v) = \vert F(u, v)\vert^2; \quad Z_{-} = F(\bar u, v)F(u, \bar v) = \vert F(\bar u, v)\vert^2.\]
Then $2$-valued group is given by
\[
x\star y = [z_1, z_2],
\]
where $z_{1,2}$ are the roots of $B(x,y,z)$. The identity is $0$ and the inverse element for $x$ is $x$ itself. Such a two-valued group is called the {\bf "square modulus"} of the original formal group.

Buchstaber classification theorem states that algebraic square modulus groups, which originate from the addition theorem for Baker–Akhiezer elliptic functions, are defined by a zero locus of the following discriminant family depending on 3 parameters 
$(a_1 ,a_2, a_3):$
\begin{equation}\label{Buniv}
B_{a_1,a_2,a_3}(x,y,z):= (x+y+z - a_2x y z)^2 -4(1+ a_3 x y z)(xy + yz+ zx + a_1 xyz).
\end{equation}
 In general, the equation $B_{a_1,a_2,a_3}(x,y,z)=0$ for such a polynomial defines the coset group $\mathcal E/\iota_{\ba}$, where $\mathcal E$ is an elliptic curve (considered as an abelian group) and $\iota_{\ba}\colon \mathcal E \to \mathcal E$ is the antipodal involution.

More precisely, consider an elliptic curve in the standard Weierstrass form $$v^2=4u^3-g_2u-g_3$$ and a point $\alpha$ on it, and set the corresponding parameters as
\begin{equation}
\label{param}
a_1=3\wp(\alpha),\quad a_2=3\wp(\alpha)^2-\frac{g_2}{4},\quad a_3=\frac{1}{4}(4\wp(\alpha)^3-g_2\wp(\alpha)-g_3),
\end{equation}
where $\wp$ is the classical Weierstrass elliptic function satisfying the equation
$$(\wp')^2=4\wp^3-g_2\wp-g_3.
$$
Then two-valued group given by the $B_{a_1,a_2,a_3}(x,y,z)$ translates into the following identity for $\wp$-functions. \cite{BuchVes}:
\begin{theorem} 
The discriminant locus multiplication law  $B_{a_1,a_2,a_3}(x,y,z) =0$ with the "elliptic parametrization" can be reduced to the addition law $X\pm Y\pm Z=0$
via the change of variables:
\begin{equation}\label{conj}
x = \frac{1}{\wp(X) +\wp(\alpha)},  y= \frac{1}{\wp(Y) +\wp(\alpha)}, z = \frac{1}{\wp(Z) +\wp(\alpha)}.
\end{equation}
\end{theorem}
This theorem follows from inspection of the original Burnside's proof for the $\wp-$function addition theorem. In the next subsection, we connect this original Burnside's proof with the coset construction corresponding to the addition law on an elliptic curve.

\subsection{$\wp-$ function, algebraic 2-valued  group law on $\mathbb P^1$ and generalization.} 
Let $\Gamma\subset \mathbb C$ be a lattice which defines an elliptic curve $\mathcal E = \mathbb C/\Gamma.$ Then the Weierstarss $\wp-$ function defines a holomorphic map 
$$
\mathbb C/\Gamma \to \mathbb P^1, 
$$
Denote the involutive automorphism of the torus $\mathcal E$ by $\tau: \tau^2 =1$.
Now we observe that  2-coset construction provides us with the projection $\pi_{1} :\mathcal E \to \mathcal E/\tau \simeq \mathbb P^1:$
$$
\pi_1(u) =\wp(u) \in \mathbb P^1.
$$
The 2-valued group operation on $\mathbb P^1$ now reads from the  2-coset construction:
if $x:=\wp(u), y:=\wp(v), \, u,v\in\mathcal E, \, x,y\in\mathbb P^1$ then 
$$
x\ast y = \pi(\pi^{-1}(x)\circ \pi^{-1}(y)),
$$ 
where $\circ-$ is the abelian group operation on the elliptic curve $\mathcal E.$
Hence we obtain
$$
x\ast y = \wp(u)\ast \wp(v),\, \pi^{-1}(x)=\pm u,\, \pi^{-1}(y)=\pm v
$$ 
$$x\ast y = \wp(\pi^{-1}(x)\circ \pi^{-1}(y)) = \wp(\pm u +\pm v)= \left[ \wp(u+v),\wp(u-v)\right]\in {\rm Sym}^2 (\mathbb C)
$$.
Any $\mathfrak S_2-$ invariant pair of numbers can be represented as a pair of roots $[z_1,z_2]$ of  a quadratic equation (which is in our case has coefficients in formal series in $x$ and $y$):
$$ 
z^2 -\Theta_1 (x,y) z + \Theta_2 (x,y) =0,
$$
\begin{equation}\label{piadd}
\Theta_1 (x,y) =\wp(u+v)+\wp(u-v) ,\, \Theta_2 (x,y) =\wp(u+v)\wp(u-v) 
\end{equation}

 \begin{theorem}\label{th:rat}
The coefficients $\Theta_{1,2}$ are rational functions in $x$ and $y$ when $x:=\wp(u), y:=\wp(v).$
 \end{theorem}
 \begin{proof}
We prove this theorem using the classical addition theorem for the Weierstrass function based on the idea of Burnside.
One can consider the beautiful Burnside determinant formula \cite{Burn} like the following system:
 \begin{equation}\label{addPi}
 {\rm det} 
 \left (\begin{array}{ccc}
 1& \wp(u)& \wp'(u)\\
 1& \wp(v)& \wp'(v)\\
 1& \wp(w)& \wp'(w)
\end{array}
\right ) = 0, \quad u+v+w =0,\, \wp'^2 = 4 \wp^3- g_2 \wp -g_3,\,  \forall u,v,w
\end{equation}
defines  the  $\wp-$function addition law. This result is a corollary of the following proposition:  
\begin{proposition}
If  $x:=\wp(u), y:=\wp(v), z:=\wp(w)$ as above then the system  \eqref{addPi} is equivalent to say that
$x,y,z$ are three roots of the cubic:
\begin{equation}\label{ABcub}
4\xi^3 -g_2 \xi -g_3  -(A\xi +B)^2 = 0,
\end{equation}
with $A:=  \frac{x' -y'}{x-y}$  , and $B = \frac{xy' -yx'}{x-y} .$
\end{proposition}
\begin{proof}
We denote $\wp(u)',\wp(v)',\wp(w)'$ by $x',y',z'$ correspondingly. Then,  the system \eqref{addPi} can be written indeed as the linear-quadratic system  with respect to $z':$

$$
\begin{cases}
(y-x)z' + (x'-y') z  = yx'-xy'\\
(z')^2 = 4z^3  - g_ 2 z -  g_3.
\end{cases}
$$
We have 

$$ 
z'= \frac{x' -y'}{x-y} z +  \frac{xy' -yx'}{x-y} = Az + B,
$$
with $A:=  \frac{x' -y'}{x-y}$ and $B:=  \frac{xy' -yx'}{x-y}.$
 We get (using the second equation of the system) that  $z$ satisfies \eqref{ABcub}. On the other hand, $x$ and $y$ are also satisfied \eqref{ABcub} which is guaranteed by the direct substitution.
\end{proof}

Then there are the following "Vi\`ete relations"
\begin{equation}\label{sigmas}
\centering
\left\{\begin{split}
\sigma_1  &=\frac{1}{4}A^2 = \frac{1}{4}\left(\frac{x' -y'}{x-y}\right)^2\\
\sigma_2 + \frac{1}{4}g_2&= \frac{1}{2}AB =  \frac{1}{2}\left(\frac{x' -y'}{x-y}\right) \left(\frac{xy' -yx'}{x-y}\right) \\ 
\sigma_3 -\frac{1}{4}g_3& = \frac{1}{4}B^2 = \frac{1}{4} \left(\frac{xy' -yx'}{x-y}\right)^2\\
\end{split}\right.
\end{equation}
and therefore

$$ 
D_{a,b,c}(x,y,z) = (\sigma_2-b)^2  - 4(\sigma_3 +c)(\sigma_1 -a) = 0
$$
 with 
 
$$b = -\frac{1}{4}g_2;\quad  c= - \frac{1}{4}g_3;\quad a =0 .
$$


By analogy with the previous example the 2-valued structure is determined as the roots of this equation

$$  
(\sigma_2+\frac{1}{4}g_2)^2 = 4(\sigma_3 - \frac{1}{4}g_3)\sigma_1
$$ 
or explicitly:

\begin{equation}
\label{lawequat}
(xy + yz + zx +\frac{1}{4}g_2)^2 = 4(xyz - \frac{1}{4}g_3)(x+y+z)
\end{equation}

We recognize immediately  the zero locus of Kontsevich polynomial $D_{0,-\frac{g_2}{4},-\frac{g_3}{4}}(x,y,z)$. But this condition still does not give us the Buchstaber formal 2-valued group law on $\mathbb C$  because taking $y=0$ we get 
$$F(x,0,z) = x^2 z^2 +(1/2 g_2 x+ g_3)z +g_3x = 0
$$instead $F(x,0,z)=(z-x)^2$.  In other hand, the locus \eqref{lawequat} gives an algebraic 2-valued group law on $\mathbb P^1$ 
when we go from the Kontsevich polynomial to the locus of the  Buchstaber one $B_{0,-\frac{g_2}{4},-\frac{g_3}{4}}(x,y,z) =0$ \cite{Buch90}.

Replacing $\mathbb C-$coordiantes $(x,y,z)$  by their $\mathbb P^1-$homogeneous counterparts $(x_1/x_0, y_1/y_0, z_1/z_0)$ and put $x_1=y_1= z_1=1$ we obtain 
in the coordinates $x_0=x, y_0 =y, z_0= z$
\begin{equation}\label{lawequat-2}
B_{0,-\frac{g_2}{4},-\frac{g_3}{4}}(x,y,z) =  (x + y+ z +\frac{1}{4}g_2 xyz)^2 - 4(1 - \frac{1}{4}g_3xyz)(xy +yz+zx) =0, 
\end{equation}
and $B_{0,-\frac{g_2}{4},-\frac{g_3}{4}}(x,0,z) = (x + z)^2 -4xz= (x-z)^2.$

This transformation changes the "unity point" $(0:1)$ for Buchstaber law to the "infinity point" $(1:0)$ for \eqref{lawequat-2}.

This law is defined by the roots $[z_1, z_2]$ of the quadratic equation \eqref{lawequat} in $z$ using \eqref{sigmas}:

\begin{equation}\label{twogroupp1}
x\ast y = \left[z_1,z_2\right] = \left[-(x +y) +\frac{1}{4}\left(\frac{x'-y'}{x-y}\right)^2, -(x +y) +\frac{1}{4}\left(\frac{x'+ y'}{x-y}\right)^2\right],
\end{equation}
where

$$\mathcal E_1: (x')^2 = 4x^3 - g_2 x -g_3, \quad \mathcal E_2: (y')^2 = 4y^3 - g_2 y -g_3
$$ describe {\it the same elliptic curve in two different
affine coordinate pairs in $\mathbb C^2: (x,x')$ and $(y,y').$} While $z-$discriminant of $D_{0,-\frac{1}{4}g_2, - \frac{1}{4}g_3}(x,y,z) $ is splitting:
$${\rm disc}_z D_{0,-\frac{1}{4}g_2, - \frac{1}{4}g_3}(x,y,z) = f(x)f(y)
$$ with $f(z) = 4z^3 -g_2 z -g_3.$

Remind that here we put $x':= \wp'(u), x:=\wp(u)$ and $y':= \wp'(v), y:=\wp(v)$ and read the Weierstrass function addition from the Burnside determinant system \eqref{addPi} as
\begin{equation}\label{WeiLaw}
\wp (u +v) = -\wp(u) - \wp(v) + \left(\frac{\wp'(u) - \wp'(v)}{2(\wp(u) - \wp(v))}\right)^2.
\end{equation}
\end{proof}

\begin{remark}
One can arrive to the Weierstrass $\wp-$ addition law almost straightforwardly starting from \eqref{twogroupp1}:

Let 
\[U:= -(x+y) +\frac{1}{4}\frac{x'^2 +y'^2}{(x-y)^2}.\]
Then 
\[z_1 = U-\frac{1}{2}\frac{x'y'}{(x-y)^2},\, z_2 = U+\frac{1}{2}\frac{x'y'}{(x-y)^2}. \]
\[z_1  + z_2 =2U = 2\Theta_1(x,y);\, z_1 z_2 = U^2 -\frac{1}{4}\frac{x'^2y'^2}{(x-y)^4} = \Theta_2(x,y),\]
and
\[z^2 -\Theta_1(x,y)z +\Theta_2(x,y) = (z-U)^2 - \frac{1}{4}\frac{x'^2y'^2}{(x-y)^4}.\]

But having in mind the relations 
$$
(x')^2 = 4x^3 - g_2 x -g_3, \, (y')^2 = 4y^3 - g_2 y -g_3.
$$
we obtain by straightforward manipulations
\[ U= \frac{1}{(x-y)^2}\left((x+y)xy -\frac{g_2}{4}(x+y) -\frac{g_3}{2}\right)\]
which implies the Buchstaber polynomial equation \eqref{lawequat-2}.
\end{remark}

In the next subsection, we discuss an extension of Burnside's determinant (also known as the {\it Frobenius-Stickelberger formula}) and motivate a generalization of the Kontsevich polynomial by considering a deformation of the $\wp$-addition theorem.

\subsection{Generalization of  algebraic 2-valued group law on $\mathbb P^1$ leads to the general Kontsevich polynomial }

In this subsection, we shall discuss a few other interesting links between the 2-valued group construction from the previous subsection and the realm of generalized Kontsevich polynomials. Technically it will lead us to a "deformation " of the addition theorem (and, hence, of the algebraic 2-group law) whose neutral element will be "shifted". It also means that the initial values in the existence and unicity theorem should be also shifted. The "shift" justifies an additional motivation to consider a more general family of Kontsevich polynomials studied in the paper.

We start with a reminder of the classical Frobenius-Stickelberger formula.  In its simplest nontrivial case, it reads as follows: if $\kappa\in \mathbb C^{*}$ a
 \begin{equation}\label{conuvw}
 u+v+w =\kappa
 \end{equation}
then,  denoting ${\tilde u}:= u-\kappa/3, {\tilde v}:= v-\kappa/3,{\tilde w}:= w-\kappa/3$ and the Frobenius-Stickelberger formula, which reads
 \begin{equation}\label{FrobStick}
 {\rm det} 
 \left (\begin{array}{ccc}
 1& \wp(u)& \wp'(u)\\
 1& \wp(v)& \wp'(v)\\
 1& \wp(w)& \wp'(w)
\end{array}
\right ) = 2\frac{\sigma(u+v+w)\sigma(u-v)\sigma(v-w)\sigma(w-u)}{\left(\sigma(u)\sigma(v)\sigma(w)\right)^3}.
\end{equation}
is equivalent to 
 \begin{equation}\label{FrobStick-2}
 {\rm det} 
 \left (\begin{array}{ccc}
 1& \wp(\tilde u)& \wp'(\tilde u)\\
 1& \wp(\tilde v)& \wp'(\tilde v)\\
 1& \wp(\tilde w)& \wp'(\tilde w)
\end{array}
\right )=  0.
\end{equation}

The condition \eqref{conuvw} crucially differs from \eqref{addPi} but  the vanishing of the Burnside determinant for  \eqref{FrobStick-2} and the Weierstrass uniformization
for \[\wp(z-\kappa/3): \wp'(\tilde u) = 4\wp(\tilde u)^2 -g_2\wp(\tilde u)-g_3\]
defines a 2-valued operation with shifted neutral element which we associate with the generalization of the Kontsevich polynomial addition law mentioned in \cite{Konts2007}.

Moreover, the appearance of Weierstrass $\sigma-$ function provides a link to the  {\it groupoid operation} which was proposed in \cite{buchstaber2005addition}. The addition law on Jacobians of genus $g\geq 1$ (hyper)elliptic curves the authors of interpret as the structure of a commutative algebraic groupoid  $\mathbb C^{3g}\to \mathbb C^{2g}.$

The Jacobian of a curve for $g=1$ can be identified with itself and the groupoid  $\mathbb C^{3}\to \mathbb C^{2}$ addition law can be interpreted as follows. Let us consider the genus 1 family of curves $y^2 = 4x^3 -g_2 x-g_3.$
Let $\mu: \mathbb C^{3}\to \mathbb C^{2}$ be the map giving by
\[\mu: (x,y,z) \to (g_2, g_3), \quad  (g_3, g_2) = (y^2 - 4x^3  + z x, -z).\]

Then we suppose that the groupoid operation for the groupoid $\mathbb C^{3}\to \mathbb C^{2}$ points $G_1 =(u,u',g_2)$ and $G_2 =(v,v',g_2)$ such that
\[\mu(G_1) = \mu(G_2) =(g_3,g_2)\] reads as
\[G_1\ast G_2 = G_3(w,w',g_2).\]
Then example 1.24 from \cite{buchstaber2005addition} provides the relation
\[w' =-\frac{uv'-u'v}{u-v} -w\frac{u' -v'}{u-v} \]
which is equivalent to vanishing of the "Burnside-Wandermonde"  determinant:
 \begin{equation}\label{group}
 {\rm det} 
 \left (\begin{array}{ccc}
 1& u& u'\\
 1& v& v'\\
 1& w& w'
\end{array}
\right )=  0.
\end{equation}

A generalization of the classical Frobenius-Stickelberger formula \eqref{FrobStick} was proposed in the paper \cite{BuchPer}. The generalization was used to solve the following functional equation in the local field of meromorphic functions:
\begin{equation}\label{BP1}
\left( f(x) + g(y) + h(z) \right)^2  = \left( F(x) + G(y) + H(z) \right), \quad x+y+z=0.
\end{equation}
Let 
\[{\tilde \mu} : \mathbb C^{3}\to \mathbb C, \quad {\tilde \mu}(u,v,w) = u+v+v = \kappa \in \mathbb C. \]
Taking $\kappa = 0$ we obtain the map 
$$\mathbb C^{3}\to \mathbb C^{2}={\tilde \mu}^{-1}(0)=\{u+v+v =0\}.$$

On the level of $\wp$-uniformization, such shift only changes the argument of $\wp$-function, but not the form of algebraic curve and preserves the groupoid structure. Such shifted groupoid structure should have a natural analogue for the genus $g$ hyperelliptic curves, which gives an extension of the results from \cite{buchstaber2005addition}. 

 A beautiful and surprising connection between  \eqref{BP1} and the algebraic addition law in the Burnside system form \eqref{addPi} confirms that $f,g,h$ solves \eqref{BP1} for given $F,G,H$ iff they satisfy the Burnside-type system
 \begin{equation}\label{fgh}
 {\rm det} 
 \left (\begin{array}{ccc}
 1& f'& f''\\
 1& g'& g''\\
 1& h'& h''
\end{array}
\right ) = 0, \quad x+y+z=0
\end{equation}
Let Weierstrass odd quasi-periodic $\zeta$-function is defined by the conditions:
$$
\begin{cases}
\frac{d\zeta(z)}{dz} = -\wp(z)\\
\lim_{z\to 0}\zeta(z) = \frac{1}{z}
\end{cases}
$$
The equation \eqref{BP1} (under some mild non-degeneracy conditions) has meromorphic solutions in $\zeta-$ functions of the form
$$\alpha\zeta(x-a_i, g_2,g_3) + \beta x + \gamma_i,\quad i=1,2,3
$$for "small letter" functions and in combinations of $\wp-$ and $\zeta-$ functions of the form
$$\alpha^2\wp(x-a_i, g_2,g_3)  + 2\gamma\alpha\zeta(x-a_i, g_2,g_3) + \frac{\gamma^2}{3},\quad i=1,2,3
$$for "big"$F,G,H$ (Main Theorem in \cite{BuchPer}).

A beautiful partial solution of \eqref{fgh} has the following form of the addition law
\begin{equation}\label{partfgh}
\left( \zeta(x) + \zeta(y) + \zeta(z) \right)^2  = \left( \wp(x) + \wp(y) +\wp(z) \right), \quad x+y+z=0.
\end{equation}

Analyzing the proof of this theorem, one can extract a "deformation" of the Weierstrass addition law \eqref{WeiLaw}
introducing (instead of $\wp$) two functions $\varphi$ and $\xi$ satisfying:
\begin{equation}\label{Deflaw}
\varphi(u+v) = \varphi(u) +\varphi(v)  + \xi(u)\xi(v) \frac{\varphi'(u) -\varphi'(v) }{\xi(u) - \xi(v)}
\end{equation}
with the initial data $\varphi(0) =\varphi'(0) = 0$ and $\xi(0) = 0$.

This deformational functional equation implies for its solution $\varphi (u)$ that the function $t(u) := \varphi'(u)$ is a solution of the following equation

\begin{equation}\label{deformKonts}
t'(u)^2 = t^3 + at^2 + bt + c, \quad t(0) = 0, t'(0) = c
\end{equation}
and 
$$a = 12\wp(\alpha),\, b=4g_2 - 48\wp^2(\alpha),\, c= 16(\wp'(\alpha))^2.
$$Here $\alpha$ is a point on the Weierstrass cubic curve with modular parameters $g_2,g_3:$
$$(\wp'(u))^2 = 4\wp(u)^3 -g_2\wp(u) -g_3
$$and the equation \eqref{deformKonts} means that 
$$t'(u)^2  = 16[4\wp(u+\alpha)^3 -g_2\wp(u+\alpha) -g_3].
$$ 
 
Taking here
 $$f(t) = t^3 +at^2 + bt +c \in \mathbb C[t],\quad a= 12\wp(\alpha),\, b=4g_2 - 48\wp^2(\alpha),\, c= 16(\wp'(\alpha))^2,$$one can write
the {\it deformed  generalized Kontsevich polynomial} $D_{a,b,c} (x,y,z)$ ( the discriminant of  $P_{a,b,c,t} (x,y,z):=f(t)-(t-x)(t-y)(t-z)$ viewed as a quadratic polynomial in $t$ is written as:
$$
D_{a,b,c} (x,y,z)={\rm disc}_t[(12\wp(\alpha) +x+y+z)t^2 - (xy + yz + xz +48\wp^2(\alpha)- 4g_2 )t + (x y z + 16(\wp'(\alpha))^2)].
$$

Again, taking $y=0$ we see that there is no "zero" or neutral element constraint $D_{a,b,c} (x,0,z)\neq (x-z)^2$ but instead of it we go to the  "inverses"  in the local chart of $\mathbb P^1$ where 
$$
x\to 1/x, y\to 1/y, z\to 1/z,\quad x\neq 0,y\neq 0, z\neq 0.
$$
We get 
\begin{multline*}
B_{a,b,c}(x,y,z) = z^2 \left((1-(4g_2 - 48\wp^2(\alpha))xy)^2 - 64(\wp'(\alpha))^2 xy(12\wp(\alpha)xy +x+y)\right)-\\
-2z\left(32(\wp'(\alpha))^2 x^2y^2 +24\wp(\alpha)xy +(4g_2 - 48\wp^2(\alpha))xy(x+y) +x+y\right) +(x-y)^2.
\end{multline*}
We see that this deformed Buchstaber polynomial defines a {\it local deformed 2-valued group law on $\mathbb P^1$} because of the  constrained condition validation:
$$
B_{a,b,c}(x,0,z) = z^2 -2zx +x^2 = (z-x)^2.
$$

\subsection{A new construction of algebraic 2-valued group law on $\mathbb C$ }
Introduce two models of elliptic curves in $\mathbb{C}^2$: 
$${\mathcal E}_1= \{(u,t) | u^2= t^3 + at^2 +bt + c\},$$
and
$$
{\mathcal E}_2= \{(u,t) | u^2=(t-x)(t-y)(t-z)\}
$$
\begin{theorem}
The following statements hold
\begin{enumerate}
\item Curves ${\mathcal E}_1$ and ${\mathcal E}_2$ has only two intersection points in the affine chart, if and only if the following equations holds
\begin{equation}\label{univlaw}
A_{a,b,c}(x,y)z^2 + B_{a,b,c}(x,y)z+ C_{a,b,c}(x,y) = 0,
 \end{equation}
where
\[ A_{a,b,c}  = (x-y)^2,\, B_{a,b,c} = -2\left[(b+xy)(x+y) +2(c+axy)\right],\,C_{a,b,c}= \left[(xy-b)^2 -4c(x+y+a)\right];\]
\item The LHS of the equation \eqref{univlaw} is $\mathfrak S_3 -$ invariant and coincides with the Kontsevich generalized polynomial \eqref{genkontspol};
\item This polynomial family $A_{a,b,c},B_{a,b,c},,C_{a,b,c},\in \mathbb C[x,y]$ defines an algebraic 2-valued mapping
\[\mu_{a,b,c}:  \mathbb C \times \mathbb C\setminus \Delta \to {\rm Sym}^2(\mathbb C) = \mathbb C^2, \mu_{a,b,c} (x,y) = [z_1,z_2]\] 
where $\Delta$ is a "diagonal" in  $\mathbb C \times \mathbb C: \Delta = \{(x,y)|x=y, x,y\in \mathbb C\}$ and
\[z_{1,2} := \frac{-B_{a,b,c}\pm \sqrt{B_{a,b,c}^2 - 4A_{a,b,c}C_{a,b,c}}}{2A_{a,b,c}}.\]
\item This mapping defines a family of 2-valued "partial"\footnote{It means that this family of laws is defined on $\mathbb C\setminus\Delta.$} group laws
\[z^2  + \Theta_1(x,y)z + \Theta_2(x,y) =0,\]
on $\mathbb C$ with rational coefficients $\Theta_{1,2}(x,y).$ 
\end{enumerate}
\end{theorem}

\begin{proof}
The proof of the first three statements is a straightforward verification. We further verify when the polynomials $A_{a,b,c}, B_{a,b,c}, C_{a,b,c}$ are simultaneously vanishing in $\mathbb C^2?$  Let $y=x$ and $A_{a,b,c}(x,x)=0$ then $B_{a,b,c}(x,x) =-2(x^3 +ax^2 +bx +c)$ and $C_{a,b,c}(x.x)=x^4 -2bx^2 -8cx + b^2 -4ac$. To find when the polynomials $B_{a,b,c}(x,x)$ and
  $C{a,b,c}(x.x)$ simultaneously vanished one should compute their resultant:
  \[ {\rm Res}[B_{a,b,c}(x), C_{a,b,c}(x)] = (a^2 b^2 - 4 b^3 - 4 a^3 c + 18 a b c - 27 c^2)^2\]
  which is the square of the first cubic polynomial discriminant $\Delta_f(a,b,c).$
  Remark that we can canonically embed $\mathbb C$ to $\mathbb P^1$ and hence
  when $\Delta_f(a,b,c)\neq 0$ the mapping family $\mu_{a,b,c}$ can be extended to a map
  \[{\tilde\mu}_{a,b,c} : \mathbb P^1\times \mathbb P^1 \to \mathbb P^2.\] 
  To verify the last assertion, we need to proceed (similar to the case of the Buchstaber polynomial in the theorem \eqref{basiclaw}).
  We obtain 
  \begin{eqnarray}\label{homunivlaw}
(XY_0 - X_0Y)^2Z^2 - \nonumber \\
- 2\left[(XY_0 +X_0Y)(bX_0Y_0 + XY)+2X_0Y_0(aXY+cX_0Y_0)\right]ZZ_0 +\\
+ \left[(XY - bX_0Y_0)^2 -4cX_0Y_0(aX_0Y_0 +XY_0 + X_0Y)\right]Z_0^2 = 0 \nonumber
 \end{eqnarray}
  \begin{itemize}
  \item  The (strong) identity element $(1:0)$ is defined.  Indeed for any $(X:X_0)\in \mathbb P^1$ the condition
  \[(-X_0)^2 Z^2 - 2XX_0ZZ_0 + X^2Z_0^2 = (ZX_0 - XZ_0)^2\]
  which means that the roots are "multiples":
  \[\mu_{a,b,c}((X:X_0),(1:0)) =  \mu_{a,b,c}((1:0),(X:X_0)) = x\ast e = e\ast x= [x,x] \]satisfies.
\item For any $(X:X_0)$ the "inverse" element is itself: indeed, this means, in our case, that
 \[\mu_{a,b,c}((X:X_0),(X:X_0)) = x\ast x,\]
 or the equation \eqref{homunivlaw} has two roots:
\[ (Z^{(1)}:Z^{(1)}_0) = (1:0),\, (Z^{(2)}:Z^{(2)}_0)= (C_{a,b,c}(x,x):2B_{a,b,c}(x,x)).\]
We see that the first solution guarantees us verification, and the second is different from $(0:0)$ when $\Delta_f(a,b,c)\neq 0.$
\item The last verification (the associativity)  is rather tedious: we need to verify that
for three com[plex numbers $x,y,z \in \mathbb C$  and for another triple $(u,v,w) \in \mathbb C$ the following
 two systems of equations:
\begin{small}
\[ 
\begin{cases}
(x-y)^2u^2 -2u \left[b(x+y) + (x+y)xy +2axy +c \right]+ (xy-b)^2 +4c(a+x+y) = 0\\
(w-z)^2u^2 -2u\left[b(z+w) + (z+w)zw +2azw +c \right]+ (wz-b)^2 +4c(a+w+z)= 0 
\end{cases}
\]
\end{small}
and
\begin{small}
\[ 
\begin{cases}
 (y-z)^2v^2 -2v \left[b(z+y) + (z+y)zy +2azy +c \right]+ (zy-b)^2 +4c(a+z+y)= 0\\
 (x-w)^2v^2 -2v \left[b(x+w) + (x+w)xw +2axw +c \right]+ (xw-b)^2 +4c(a+x+w)= 0.
\end{cases}
\]
\end{small}
after the elimination of $u$ and $v$ define the same 4-set $(x, y, z,w)\in \mathbb C^4.$
\end{itemize}
 \end{proof}
 Here we introduced construction of a family of 2-valued algebraic mappings
 \[\mu_{\bar a}: \mathbb C\times \mathbb C \to {\rm Sym}^2(\mathbb C) = \mathbb C^2, \]
 parametrized by $\bar a =(a,b,c)\in \mathbb C^3.$

\subsection{Algebraic coset 2-valued group structures on $\mathbb P^1$}\label{subs:alg}

We restrict our attention to {\it algebraic} 2-valued groups with the multiplication law given by the equation $F(x,y,z)=0,$ where $F(x,y,z)$ is a symmetric polynomial which has degree 2 in each
variable and such that $F(x,y,0)= (x-y)^2$  All such groups are commutative and involutive.
\begin{remark}
The condition $F(0,y,z)= (z-y)^2$ which follows from the symmetry condition means that $0$ is the "strong group unity".
\end{remark}
  
    This section is devoted to a new description of the coset  2-valued group structure on the projective line $\mathbb P^1.$ Such a structure was obtained for the first time in the paper \cite{BuchRees} (see also \cite{Buch06}) by the "square moduli" construction from the formal group whose exponent is the Jacobi  ${\rm sn-}$ elliptic sine function.
   
 The first author and A. Veselov \cite{BuchVes} had shown that the equation $F(x,y,z)=0$ with symmetric polynomial with  $F(x,y,0)= (x-y)^2$  defines a 2-valued group 
  iff the polynomial $F(x,y,z)$ coincides with the Buchstaber polynomial \eqref{buchstpol1}. We should stress that for {\it generic} values of the parameters $a_1,a_2,a_3$ all
 their construction defines 2-valued formal groups on $\mathbb C$.   Here we shall interested in the natural question: when such 2-valued formal group structures could be extended
 to $\mathbb P^1$? The answer is quite natural and nice – it is possible if and only if the "uniformizing" elliptic curve $\mathcal E_{a_1,a_2,a_3}: v^2 = u^3 + a_1u^2 + a_2 u + a_3$
 is a non-degenerate.

Our aim here is to propose a 2-valued group law description (similar to \cite{BuchRees} and  \cite{Buch06}) based on an algebraic mapping $\mathbb P^1\times \mathbb P^1 \to \mathbb P^2$ starting with the Buchstaber polynomial family and to find algebraic conditions on the parameters $a_i, i=1,2,3$ which define various famous examples of algebraic 2-valued groups. 
The formal group law behind this construction has as the exponent the function
\[\frac{\sigma(u)\sigma(\alpha)}{\sigma(u+\alpha)}\exp(-\zeta(\alpha)u),\]
where $\sigma(u)$ is the Weierstrass $\sigma-function$ and $\zeta(z) =\frac{d\log(\sigma(z)))}{dz}.$

Consider B-polynomial in the form
\begin{equation}\label{eq:Bpol47}
B_{a_1,a_2,a_3}(x,y,z) = (x+y+z - a_2x y z)^2 -4(1+ a_3 x y z)(xy + yz+ zx + a_1 xyz).
\end{equation}
\begin{theorem}\label{basiclaw}
The following statements hold:
\begin{enumerate}
    \item polynomial \eqref{eq:Bpol47} defines 2-valued algebraic group structure on $\mathbb{C}$ for any choice of $(a_1,a_2,a_3)$,
    \item this 2-valued algebraic group structure extends to $\mathbb{P}^1$ iff  $t$-discriminant of the cubic polynomial  \[ f(t) = t^3 + a_1t^2 + a_2 t + a_3\]
    is not zero.
\end{enumerate}
\end{theorem}
\begin{proof}
To prove the first statement in general case it is enough to verify associativity conditions over the open dense subset in $\mathbb{C}^3$ with coordinates $a_1,a_2$ and $a_3$.  In work \cite{Buch90} it was shown that if we take this open dense subset to be a complement to the zero loci of discriminant for $f(t)$, then associativity condition is satisfied. Since associativity condition translates into algebraic problem of coincidence of multi-valued polynomials in $a_1,a_2$ and $a_3$, this equality extends to the whole $\mathbb{C}^3$ (see equations \eqref{eq:11}, \eqref{eq:12} and proposition \ref{prop:3.3}).

Now let us switch to the second statement. Let $(x_1: x_0)$ be $\mathbb P^1-$homogeneous coordinates such that $|x_1|^2 + |x_0|^2 \neq 0, x_k \in \mathbb C, k=0,1.$
Consider a family of lines for a pair $(z_0,z_1)\in \mathbb C^2$
\begin{equation}\label{linfam}
x_0 z_1 + x_1z_0 =0.
\end{equation}
We choose homogeneous coordinates $(u_2:u_1:u_0)$  on the projective. plane $\mathbb P^2$ 
such that $u_k \in \mathbb C, k=0,1,2$ and $|u_2|^2 +|u_1|^2 + |u_0|^2 \neq 0.$

It is well–known that  the symmetric square ${\rm S}^2(\mathbb P^1)$ can be identified with $\mathbb P^2$ by the multiplication 
\[\mu: \mathbb P^1 \times \mathbb P^1 \to (\mathbb P^1 \times \mathbb P^1)/\mathfrak S_2 \simeq {\rm S}^2(\mathbb P^1)  = \mathbb P^2. \]

The image of the linear family pair
\[ (x_0 z_1 + x_1z_0 , y_0 z_1 + y_1z_0 )\]
(see  \eqref{linfam}) in $\mathbb P^1$ to the family of quadrics in $\mathbb P^2:$
\begin{equation}\label{lsqfam}
u_2 z_0^2+ u_1z_0z_1 +u_0 z_1^2 =0,
\end{equation}
where
\[u_2 = x_1y_1; u_1 = x_1y_0 +x_0y_1, u_0 = x_0y_0.\]
 
We provide the multiplication law $\mathbb P^1\times \mathbb P^1 \to \mathbb P^2$ with the Buchstaber polynomial in the following manner: we make the change of variables in the 
"affine" Buchstaber polynomial using the homogeneous coordinates in $\mathbb P^1.$
\[ x= \frac{x_1}{x_0}, y= \frac{y_1}{y_0}, z= \frac{z_1}{z_0}, x_0\neq 0,y_0\neq 0, z_0\neq 0.\] and
\begin{multline*}
B_{a_1,a_2,a_3}\left(\frac{x_1}{x_0},\frac{y_1}{y_0},\frac{z_1}{z_0}\right) = \left(\frac{x_1}{x_0}+\frac{y_1}{y_0}+\frac{z_1}{z_0} - a_2\frac{x_1}{x_0}\frac{y_1}{y_0}\frac{z_1}{z_0}\right)^2 - \\
-4\left(1+ a_3 \frac{x_1}{x_0}\frac{y_1}{y_0}\frac{z_1}{z_0} \right)\left(\frac{x_1}{x_0}\frac{y_1}{y_0} + \frac{y_1}{y_0}\frac{z_1}{z_0}+ \frac{x_1}{x_0}\frac{z_1}{z_0} + a_1 \frac{x_1}{x_0}\frac{y_1}{y_0}\frac{z_1}{z_0}\right).
\end{multline*} 

To verify that the image of this map falls in $\mathbb P^2$ we need to check the following 
 system of compatibility conditions for the map
\[(x_0 z_1 + x_1z_0 , y_0 z_1 + y_1z_0 ) \to u_2 z_0^2+ u_1z_0z_1 +u_0 z_1^2 :\]
$$
\begin{cases}\label{compsys}
u_2 = (x_1y_0 - x_0y_1)^2 = 0;\\
u_1 = a_2 x_1^2 y_1 y_0  +  a_2 y_1^2 x_1 x_0 - y_0^2 x_1 x_0   - x_0^2 y_1y_0   - 2a_1 x_1x_0 y_1 y_0  -2a_3 x_1^2 y_1^2 = 0;\\
u_0 = \Omega(x_1y_1)^2 +(x_0y_0)^2  - 2a_2 x_0x_1y_0y_1- 4a_3 x_0x_1y_1^2 - 4a_3 x_1^2y_0y_1 =0,\\
\end{cases}
$$
where $\Omega : = a_2^2 - 4a_1a_3$. The system \eqref{compsys} is incompatible if
\begin{itemize}
\item $256((a_1^2 - 4a_2\Omega + 2a_1a_2a_3 - 27a_3^2)^2 \neq 0.$
\item The last expression (up to the numeric factor) coincides with $(\Delta_f (a_1,a_2,a_3))^2$ (see \cite{GKZ}, ch.12  (1.34)).
\end{itemize}

The map $\mu$  is not well defined at the point $ ((x_0 : x_1), (y_0 : y_1))$ if
all coordinates $u_k, k=0, 1, 2$ are equal to $0.$
\begin{enumerate}
\item The vanishing of the first is equivalent to the condition $x_1y_0 = x_0y_1.$
Substituting $x_0y_1$ instead of $x_1y_0$  in the second and in third equations we obtain
$$
\begin{cases}
\Omega(x_1y_1)^2 +(x_0y_0)^2  = 2a_2 (x_0y_1)^2 + 8a_3 x_0x_1y_1^2;\\
y_1(-a_2 x_1 x_0  + a_3 x_1^2  + a_1  x_0^2 ) =x_0y_0.
\end{cases}
$$
\item Assuming $x_0y_0\neq 0$ and denoting $k:= \frac{x_1}{x_0} = \frac{y_1}{y_0}$ we obtain:
\begin{eqnarray}\label{twopols}
\Omega k^4 + 1  = 2a_2 k^2 + 8a_3 k^3;\\
k(-a_2 k  + a_3 k^2  + a_1) = 1.
\end{eqnarray}
To check when these quartic and cubic polynomials have a common zero locus one can compute their resultant:
\begin{equation}\label{resultant}
R(a_1,a_2,a_3)  = 256((a_1^2 - 4a_2)\Omega + 2a_1a_2a_3 -27a_3^2)^2 =  256 (\Delta_f(a_1,a_2,a_3))^2 ,
\end{equation}

 The condition $R(a_1,a_2,a_3) \neq 0$ guaranties that there are no common zeros of two polynomials \eqref{twopols} and the image of the multiplication 
 falls into $\mathbb P^2.$  It defines a 2-valued group structure on $\mathbb P^1$ if and only if
the elliptic curve $\mathcal E_f : v^2 = u^3 + a_1u^2 + a_2 u + a_3$  is non-degenerate.
 We observe that if $\Omega \neq 0$ the choice $a_3 =0$ implies still non-degenerate cubic if $R(a_1,a_2,0) = (a_1 a_2)^2 - 4a_2^3 \neq 0.$ 
  It is interesting to consider the partial case of our condition when $a_1=0.$ In this case the condition transforms in a non-triviality of the 
  cubic discriminant for the Weierstrass elliptic curve $v^2 = u^3 + a_2 u + a_3:$
\[R(0,a_2,a_3) =  4a_2^3 + 27a_3^2 \neq 0.\] 
If as above $a_3 =0$ we still have the non-trivial non-degeneracy condition $R(0,a_2,0) =  4a_2^3 \neq 0.$ 
\item Another partial case of a similar reduction to a Weierstrass cubic appears when 
 $\Omega = 0.$ The resultant \eqref{resultant} takes the form
\[R(a_1,a_2,a_3)_{\Omega=0} = 256\left(\frac{a_2^3}{2} - 27a_3^2\right)^2.\]
and implies the cubic discriminant condition  $\frac{a_2^3 }{2} - 27a_3^2\neq 0 $ which is almost identical (up to rescaling) to the condition $R(0,a_2,a_3) =  0.$

Remark further that the constraint $a_3 =0$ in the latter case is impossible because it implies $a_2 =0$  which contradicts the first equation of \eqref{twopols}.
\item Let $x_0y_0$ be equal to $0.$ Then it is easy to conclude that this condition implies that both coordinates should equal to $0$.
 \end{enumerate}
\end{proof}

This result is a full analogy of Th.12, p.19 \cite{Buch06} and prop. 41, p.334 \cite{BuchRees}. While the (non-degenerate) elliptic curve behind the 
Jacobi sn-function uniformizes coset 2-valued group construction in these papers, the elliptic curves,  which are related to the similar construction
in our Theorem, depend on the points of the quartic surface $R(a_1,a_2,a_3) = 0$ in the frame of the "universal parameter uniformization" \eqref{unipar}.

The group laws for degenerated curves (nodal $v^2 = u^3 + a_1u^2,  a_1\neq 0, a_2=a_3=0$ and cuspidal $v^2 = u^3, a_k =0, k=1,2,3$) defines 2-valued formal 
 group laws only on $\mathbb C$ which are not extendable to $\mathbb P^1.$

The nodal cubic group 2-valued group law exactly corresponds to Example 2.3 and the Buchstaber polynomial $B_{(a_1,0,0)}$ after an appropriate change of coordinates and 
the normalization choice $a_1 =1$ coincides with the Mordell modification \eqref{mord} of generalized Markov cubic polynomial and for arbitrary $a_1$ the formal group law has 
appeared in the algebraic topology in the context of K-theory. The (square root of) this Buchstaber polynomial is exactly the denominator in the product law for degree 0 Gegenbauer polynomials:

\[ 
G_0 (x) G_0(y) = \frac{1}{2\sqrt \pi} \int_{-1}^{1}\frac{G_0 (z) dz}{\sqrt{(1-x^2 - y^2 - z^2 + 2xyz)}}. 
\]


The cuspidal curve case with $a_1=a_2=a_3=0$ coincides with the group law from Example 2.2. This law also can be found in the context of the usual cohomology theory.
This case is quite interesting from the {\it multiplication kernels} perspectives mentioned in the Introduction. The (square root of) the Buchstaber polynomial is exactly the 
the denominator in the product law for degree 0 modified Bessel functions of the first kind (the famous "Sonine–Gegenbauer formula"):

\begin{multline}
I_0 (2\sqrt x) I_0(2\sqrt y) = \frac{1}{\pi} \int\limits_{x -2\sqrt{xy} +y}^{x+2\sqrt{xy}+y}\frac{I_0 (2\sqrt z) dz}{\sqrt{-B_{(0,0,0)}(x,y,z)}} =\\
= \frac{1}{\pi} \int\limits_{x -2\sqrt{xy} +y}^{x+2\sqrt{yx}+y}\frac{I_0 (2\sqrt z) dz}{\sqrt{2(xy +xz+yz) -(x^2+y^2+z^2)}} ,\ x,y>0.
\end{multline}


The classification of all algebraic 2-valued group structures on $\mathbb P^1$ for the first time was obtained in \cite{Kluwer}. In this paper, it was also shown (using the result from \cite{Buch2}) that the 2-valued group structure on $\mathbb P^1$ is determined by three parameters of a {\it cubic polynomial} ${\tilde f} = 2+{\tilde a_1}x+{\tilde a_2}x^2 +{\tilde a_3}x^3 $ (ch. 4 of  \cite{Kluwer}).

\subsection{M\"{o}bius group action and Buchstaber-Kontsevich polynomials}\label{subs:moeb}
 
Recall that a {\it M\"{o}bius Transformations}  is a rational function of degree one, so that as a transformation $F$ of the extended complex plane $\bar{\mathbb C}$
\[z \to \gamma(z) =\frac{Az +B}{Cz+D}, 
\left(\begin{matrix}
A & B \\
C & D
\end{matrix}\right)\in SL_2(\mathbb C), AD - BC =1; A,B,C,D\in \mathbb C.
\]
Here $F(\infty) = \frac{A}{C}$ and $F(-\frac{D}{C}) =\infty.$ 

The group of M\"{o}bius  Transformations is exactly thrice transitive:
if  $u,v,w$ are three points of $\mathbb C$, then
\[
z \to \frac{w-v}{w-u}\frac{z-u}{z-v}
\]
sends $u$ to $0$, $v$ –  to $\infty$ and $w$ to $1.$

On the other hand, by direct computation, it is not hard to show
that the only M\"{o}bius Transformation which fixes $( 0, 1,\infty)$  is the $\gamma \equiv \mathbb Id.$
The following lemma is "common knowledge".
\begin{lemma}
The group of M\"{o}bius transformations that preserves $\{0,1,\infty \}$ is isomorphic to $\mathfrak S_3$, generated by 
\[z\to \frac{1}{z},\quad z\to 1-z.\]
\end{lemma}

We shall write the M\"{o}bius action to $\mathbb P^1$ as 
\[\gamma: \mathbb P^1 \to \mathbb P^1 \] as
\[
[X:X_0] \to [ AX+BX_0 : CX+DX_0].
\]
Take three points in  $\mathbb P^1: ,\ z_1 = (0:1), z_2= (1:1), z_3 = (1:0).$  The matrix family
\[
\gamma_1(C) =
\left(\begin{matrix}
1 & 0 \\
C & 1
\end{matrix}\right)
\]
form a  1-parametric subgroup $G_1\subset SL_2(\mathbb C)$ and $\gamma_1 (C) (0:1) = (0:1)$ which means that $G_1-$ the stationary subgroup of the point $z_1$.
Similarily, $G_2\subset SL_2(\mathbb C)$ such that
\[
\gamma_2(A) =
\left(\begin{matrix}
A & 1-A \\
A-1 & 2-A
\end{matrix}\right)
\]
form a  1-parametric stationary subgroup for the point 
\[
z_2: \gamma_2 (A)(z_2) = (A + 1-A: A-1 +2-A)= (1:1)
\]
and 
$G_3\subset SL_2(\mathbb C)$ such that
\[
\gamma_3(B) =
\left(\begin{matrix}
1& B \\
0 & 1
\end{matrix}\right)
\]
form a  1-parametric stationary subgroup for the point $z_3: \gamma_3 (B)(z_3) = (1: 0).$
One can easily extend the diffeomorfism action $\gamma : \mathbb P^1 \to \mathbb P^1$ obtaining a new multiplication map
$\mu^{\gamma} : \mathbb P^1 \times \mathbb P^1 \to \mathbb P^2$ and to check the commuttativity of the diagram

$$\xymatrix{\mathbb P^1\times\mathbb P^1\displaystyle{\ar[r]^{\mu}}\displaystyle{\ar[d]_{\gamma\times \gamma}}&
\mathbb P^2\displaystyle{\ar[d]^{\gamma}}\\
\mathbb P^1\times\mathbb P^1 \displaystyle{\ar[r]^{\mu^\gamma}}&\mathbb P^2.}$$
Now it is easy to study how the elementary 2-valued group structure in $\mathbb C$ can be "extended" to 2-valued group structures in $\mathbb P^1:$
\begin{proposition}
Take the triple $(0,1,\infty)$ of points in $\mathbb P^1$.
\begin{enumerate}
\item  "neutral element" is $\infty = (1:0),$
the orbit is $SL_2(\mathbb C)/G_3$ and corresponds to Kontsevich polynomial;
\item  "neutral element" is $0 = (0:1),$ 
the orbit is $SL_2(\mathbb C)/G_1$ and corresponds to Buchstaber polynomial;
\item "neutral element" $(1:1),$
the orbit is $SL_2(\mathbb C)/G_2$  and corresponds to the polynomial $P(A)$
such that
\[P(A) =\gamma_B^{*}(B_{a_1,a_2,a_3}) = \gamma_D^{*}(D_{a,b,c}),  \]
with
\[
\gamma_B=
\left(\begin{matrix}
A& 1 \\
A-1 & 1
\end{matrix}\right),
\gamma_D =
\left(\begin{matrix}
1 & B \\
1 & B+1
\end{matrix}\right)
\]
\end{enumerate}
\end{proposition}
\begin{cor}
Three extended multiplication laws on $\mathbb P^1$ are described as follows:
\begin{enumerate}
\item The multiplication law $\mu_{(a,b,c)}^{1}$ is defined by the Buchstaber polynomials, which correspond to the M\"{o}bius group action orbit $\mathcal O_{(a,b,c)}^{1} = SL_2(\mathbb C)/G_1$
of the strong neutral element $(0:1);$
\item  The multiplication law $\mu_{(a,b,c)}^{2}$ is defined by the Kontsevich polynomials, which correspond to the M\"{o}bius group action orbit $\mathcal O_{(a,b,c)}^{3} = SL_2(\mathbb C)/G_3$
of the strong neutral element $(1:0);$
\item  The multiplication law $\mu_{(a,b,c)}^{3}$ is defined by the polynomial $P(A)$, which correspond to the M\"{o}bius group action orbit $\mathcal O_{(a,b,c)}^{2} = SL_2(\mathbb C)/G_2$
of the strong neutral element $(1:1).$
\end{enumerate}
\end{cor}

\section{Correspondence between Buchstaber and Kontsevich  polynomials}\label{sec:BK}
In this section, we give a geometric description of the correspondence between Buchstaber and Kontsevich polynomials. These polynomials carry a lot of similar properties and even coincide for some very specific cases like a kernel for some degeneration of $D2$ equation\footnote{The general definition of $DN–$ operators and corresponding equations see \cite{GolSt}. In particular, a generic regular $D2$ equation 
\[ L\phi_{\lambda}(t) = \lambda \phi_{\lambda}(t)\]  (up to some elementary transformations) is given by the following Heun–type operator $L  = f (t)\frac{d^2}{dt^2} + f(t)\frac{d}{dt} + t$, where $f(t) = t^3 + At^2 + Bt.$ The parameter $\lambda$ plays the role of an accessory one.} 
In fact the following theorem holds
\begin{theorem}\label{th:KB1}
Generalized Kontsevich polynomial corresponds to Buchstaber polynomial with parameters $a_1=a,\,a_2=b,\,a_3=c$ after the following change of variables
\begin{equation}
    (xyz)^2D(a,b,c|-x^{-1},-y^{-1},-z^{-1}) = B_{a,b,c}(x,y,z).
\end{equation}
This means, in particular, the Buchstaber polynomial becomes a discriminant type polynomial of the form
\begin{equation}\label{eq:KB1}
    B_{a,b,c}(x,y,z) = {\rm disc}_t\left[xyz(t^3+at^2+bt+c)-xyz(t+1/x)(t+1/y)(t+1/z)\right]
\end{equation}
\end{theorem}
\begin{proof}
    Direct computation. We see that both polynomial families transform one to another on an open everywhere dense set.
\end{proof}
\begin{remark}
    The discriminant formula \eqref{eq:disc1} represents a discriminant of Lagrange polynomial for cubic polynomial $P(t) = t^3+at^2+bt+c$ crossing the points
    \[
    (x,P(x)), \quad (y,P(y)), \quad (z,P(z)).\] Indeed, since we choose $3$ points Lagrange interpolation polynomial should have degree $2$, which is exactly the degree of 
    $$    \phi(t) =t^3+at^2+bt+c-(t-x)(t-y)(t-z).
    $$    On the other hand, one checks that $\phi(x)=P(x)$, $\phi(y)=P(y)$ and $\phi(z)=P(z)$. Setting the discriminant of this polynomial to $0$ is equivalent to the constraint that the corresponding points on the elliptic curve belong to the line, which is exactly the geometric form of Abel's law on the elliptic curve
    
\end{remark}
Both Buchstaber and Kontsevich and polynomials enjoy a split discriminant property, which was mentioned before. Such property plays important role in the theory of integrable systems and was well studied in \cite{DK}:

\begin{theorem}
The generalized Kontsevich polynomials $D_{a,b,c}(x,y,z)$ are
{\it strongly discriminantly separated} (V. Dragovi\'c, K. Kuki\'c \cite{DK})
which means that there exists a polynomial $f = f(u)$ such that the discriminant of the polynomial $D_{a,b,c}(x,y,z)$ in z as a polynomial in variables $x$ and $y$ factorizes:
$${\rm disc}_z (D_{a,b,c})(x, y) = f(x)f(y)
$$
and moreover
$${\rm disc}_x (D_{a,b,c})(y, z) = f(y)f(z)
$$
and
$${\rm disc}_y (D_{a,b,c})(z, x) = f(z)f(x)
$$
\end{theorem}

 In contrast to this study, where the authors were mostly interested in the moduli spaces of elliptic curves (i.e. families parametrized by modular parameters), here we consider unreduced families and the following statement holds
\begin{theorem}
All symmetric polynomials in 3 variables $S(x,y,z)$ of degree $2$ in each variable, which discriminant splits
\begin{equation}
    {\rm disc}_z(S(x,y,z)) = 16f(x)f(y),\quad 
\end{equation}
and it satisfies the formal 2-group identity property, i.e.
\begin{equation}
    S(x,0,z)=(x-z)^2.,
\end{equation}
coincide with Buchstaber-Kontsevich-type polynomials
\begin{equation}
    (x+y+z - a_2x y z)^2 -4(1+ a_3 x y z)(xy + yz+ zx + a_1 xyz).
\end{equation}
\end{theorem}
Such correspondence provides a deep relation between integrability, formal 2-group theory, and the theory of generalized shift operators. Moreover, such a connection brings us to an involution on the space of corresponding elliptic curves. In the following subsection, we give a geometric interpretation of the Kontsevich polynomial and study the geometry of this involution.

\subsection{Intersection of elliptic curves and Buchstaber-Kontsevich involution}
The discriminant formula for the generalized Kontsevich polynomial provides an idea of the geometric problem behind this computation. Obviously, the loci itself connected with Abel's law on an elliptic curve. However, one may consider it as an elementary problem in the intersection theory. 

Consider two families of elliptic curves $u^2=t^3+P_2(t)$. The first family $\mathcal{W}$ is parametrized by the coefficients of the quadratic part $P_2(t)$ and reads as
\begin{equation}
    \mathcal{C}_{a,b,c}:\quad u^2=t^3+at^2+bt+c.
\end{equation}
The second family $\mathcal{L}$ is parametrized by the roots of $t^3+P_2(t)$ and reads
\begin{equation}
    \mathcal{E}_{x,y,z}:\quad u^2=(t-x)(t-y)(t-z).
\end{equation}
These two families correspond to the "unreduced" models of the standard moduli spaces $\mathbb{H}/SL_2(\mathbb{Z})$ and $\mathcal{M}_{\mathcal{E}}$. That means that the quotient with respect to the M\"obius action of the complement to the discriminant divisors 
$$4 a^{3} c -a^{2} b^{2}-18 abc +4 b^{3}+27 c^{2}=0 \quad \text{and}\quad
(x - z)^2(x - y)^2(y - z)^2=0
$$are isomorphic to $\mathbb{H}/SL_2(\mathbb{Z})$ and $\mathcal{M}_{\mathcal{E}}$ correspondingly.

For a generic choice of $a,b,c$ and $x,y,z$ elliptic curves $\mathcal{C}_{a,b,c}$ and $\mathcal{E}_{x,y,z}$ intersect in $4$ points in the chosen affine chart. From this point of view, generalized Kontsevich polynomial 
\begin{equation}\label{eq:disc}
    D(a,b,c|x,y,z) = \operatorname{disc}_t\left[t^3+at^2+bt+c-(t-x)(t-y)(t-z)\right]
\end{equation}
gives a locus in $\mathcal{W}\times \mathcal{L}$ where the number of intersection points drops from $4$ to $2$.

Theorem \ref{th:KB1} states that Kontsevich polynomial transfers to Buchstaber one by switching from roots of elliptic curve to their inverses multiplied by $-1$, i.e.
\begin{equation}\label{invol}
    \sigma:\quad x\rightarrow-\frac{1}{x},\quad y\rightarrow-\frac{1}{y},\quad z\rightarrow-\frac{1}{z}
\end{equation}
in the chart when $x\neq 0,y\neq 0, z\neq 0.$
Such change of coordinates defines an involution on $\mathcal{M}_{\mathcal{E}}$ which divides discriminant by $(xyz)^2$.

\begin{theorem}
Fixed points of the involution \eqref{invol} form a subfamily in $\mathcal{L}$ which is given by the union of rational curves in $\mathbb{P}^2$. It equation reads as
\begin{equation*}
    \left(x z -y^{2}\right) \left(x y -z^{2}\right) \left(x^{2}-y z \right) \left(x^{2} y  +z^{2} x +y^{2} z -3 x  yz\right) \left(x^{2} z +y^{2} x +y \,z^{2}-3 x  yz\right)=0
\end{equation*}
where cubic factors 
$$x^{2} y  +z^{2} x +y^{2} z -3 x  yz=0,\quad x^{2} z +y^{2} x +y \,z^{2}-3 x  yz=0
$$are covered by singular fibers of the Hesse cubic family
$$v^3+w^3+r^3-3 \chi_i vwr=0, \quad
$$where $\chi$ is a cubic root of unity.
\end{theorem}
\begin{proof}
The explicit form of the invariant loci may be found straightforwardly in the equation
\begin{equation}
    j(x,y,z)=j(-1/x,-1/y,-1/z),
\end{equation}
where $j(x,y,z)$ is a $j$-invariant of the curve $\mathcal{E}_{x,y,z}$. Now let us show that cubic terms are indeed connected with the Hesse cubic. Using the change of coordinates 
\begin{equation}
    x^2y=v^3,\quad y^2z=w^3\quad z^2x=r^3.
\end{equation}
This is not birational transformation, however, the factors of invariant loci transform to the rational functions over $\mathbb{Z}[\exp(2\pi i/3)]$:
\begin{equation}
    x^{2} y  +z^{2} x +y^{2} z -3 x  yz = v^3+w^3+r^3-3 \chi_i vwr,\quad \chi_i^3=1,\quad i=1..3
\end{equation}
\begin{align*}
    x^{2} z +y^{2} x +y \,z^{2}-3 x  yz=\\= \frac{\left(p u +p v +u v \right) \left(u^{2} p^{2}-p^{2} u v +p^{2} v^{2}-p \,u^{2} v -p u \,v^{2}+u^{2} v^{2}\right)}{v u p}.
\end{align*}
\begin{align*}
    \left(x z -y^{2}\right) \left(x y -z^{2}\right) \left(x^{2}-y z \right)=\\=-\frac{\left(u -v \right) \left(u^{2}+u v +v^{2}\right) \left(p -v \right) \left(p^{2}+p v +v^{2}\right) \left(p -u \right) \left(p^{2}+p u +u^{2}\right)}{v u p}
\end{align*}
One can see a Hesse-type cubic appears after the change of coordinates. One can check that this cubic is singular. Singularities up to the action of the permutation group are given by
\begin{equation}
    [\zeta_j^2\xi_i:\zeta_j\xi_i:1],
\end{equation}
where $\zeta_j$ and $\xi_i$ are the independent cubic roots of identity. 

On the other hand, one may straightforwardly show that both curves 
\begin{equation}
    x^{2} y  +z^{2} x +y^{2} z -3 x  yz=0\quad \text{and} \quad x^{2} z +y^{2} x +z^{2}y -3 x  yz=0
\end{equation}
are rational projective curves in $\mathbb P^2$. Indeed, both curves have a singular point at $[1:1:1]$. This allows to construct the rational parametrization. For the first curve 
$$x^{2} y  +z^{2} x +y^{2} z -3 x  yz=0,
$$such parametrization reads as
\begin{equation}
    x = \left(u +v \right)^{2} v,\quad  y = -u^{2} \left(u +v \right), \quad z =
v^{2} u.
\end{equation}
For the second curve, the steps are the same since it transforms to the first one under cyclic permutation
$$x\to z,\quad  y\to y,\quad  z\to x,
$$
which finishes the proof.
\end{proof}

\begin{remark}{\rm {\bf Elliptic modular surfaces from Hesse-like cubics.}}

We have also in mind the construction of {\it elliptic rational surfaces} in $\mathbb P^3$ related to our Hesse-type cubics. Their description (see i.e.\cite{iskovskikh1996algebraic} and \cite{Meyer}) is the follows:

   Let $\Gamma\subset {\mathbb PSL}_2 (\mathbb Z)$ be a torsion-free congruence 
   subgroup of order bigger than 2 and let $ \Sigma_{\Gamma}:=(\mathbb H/\Gamma)^*$ be 
 the corresponding modular curve. Than we shall denote by $\pi: S_{\Gamma}:\to \Sigma_{\Gamma}$ 
 the universal family of elliptic curves which is usually called {the elliptic modular surface associated with $\Gamma.$ }
The appeared invariant rational surfaces can be obtained (see \cite{beauville1996complex}) as a blow-up of some surface $\mathcal S_{\Gamma}\subset \mathbb P^2\times\mathbb P^1.$
The surface $\mathcal S_{\Gamma}$ is an elliptic fibration over $\mathbb P^1$  having four singular fibers:
\begin{enumerate}
\item The Hesse pencil 
\[\mathcal S_{\Gamma}: x^{3}   +z^{3}  +y^{3}  -3\lambda x  yz=0\quad \text{with}\quad \Gamma = \Gamma(3)\] and four fibers of type $I_3$ in the Beauville classification;
\item  Our cubic pencil
\[\mathcal S_{\Gamma}: x^{2}y   +z^{2}x  +y^{2}z  -3\mu x  yz=0\quad \text{with}\quad  \Gamma = \Gamma_0(9)\cap \Gamma_1(3)\] and three fibers of type $I_3$ and one of the type
$I_9$.
\end{enumerate}

\end{remark}

\subsection{ $\star$-involution and generalized Buchstaber-Kontsevich polynomials}\label{subs:var}
 

 Since both Buchstaber and Kontsevich polynomial are symmetric in $x,y$ and $z$, one may rewrite via generators of the symmetric polynomial ring. We start with the general case when the number of the polynomial ring generators is $n\in \mathbb N:$
$$ \mathbb C[x_1,\ldots, x_n]^{\mathfrak{S}_n} \simeq \mathbb C[\sigma_1, \ldots, \sigma_n], \quad \sigma_1 = x_1+\ldots+x_n,\ldots \sigma_n = x_1...x_n.
$$
In order to allow Buchstaber-Kontsevich involution, one needs to extend the polynomial ring $\mathbb C[\sigma_1, \ldots, \sigma_n]$ to the ring of Laurent polynomials which is simply localisation with respect to $\sigma_n$, i.e. $\mathbb C[\sigma_1, \ldots, \sigma_{n-1},\sigma_n^{\pm 1}]$.
Define the map 
 
$$\star : \mathbb C[\sigma_1, \ldots, \sigma_{n-1},\sigma_n^{\pm 1}]\to \mathbb C[\sigma_1, \ldots, \sigma_{n-1},\sigma_n^{\pm 1}]
$$ 
 by the formula
 
$$ \mathcal P \to \mathcal P^{\star}, \quad \mathcal P^{\star}(\sigma_1,\ldots,\sigma_n) := \mathcal P\left(\sigma_1/\sigma_n,\ldots, \sigma_{n-1}/\sigma_n,1/\sigma_n\right).
$$
The following lemma holds
 \begin{lemma}\label{star}
 The $\star-$map induces an involution on the invariant Laurent polynomials due to the isomorphism of rings
$$\mathbb C[x_1^{\pm 1}, \ldots, x_n^{\pm 1}]^{\mathfrak S_{n}} \simeq \mathbb C[\sigma_1, \ldots,\sigma_{n-1}, \sigma_n^{\pm 1}].
$$
 \end{lemma}
 Though the case of general $n$ has also some interesting applications (for example when we study polynomial Poisson structures in the frame of "$\star-$algrbras) we shall 
 restrict ourselves in this paper to the case of $n=3$ ( most relevant to the families of Buchstaber–Kontsevich polynomials).

 Recently, G. Cotti and A. Varchenko \cite{cotti2020markov} have observed that for $n=3$ the involution (\ref{star}) provides a "$\star-$ generalization of the famous Markov cubic relation
 
$$ p^2 + q^2 +r^2 = 3pqr, \quad p,q,r\in \mathbb Z
$$
 by the following polynomial identity ( they call it "$\star-$Markov Laurent polynomial equation"):
 
$$
pp^{\star} + qq^{\star} +rr^{\star} -pqr = \frac{3\sigma_1\sigma_2 - \sigma_1^3}{\sigma_3},\quad p,q,r\in \mathbb Z[\sigma_1,\sigma_2,\sigma_3^{\pm 1}].
$$
  This generalization was motivated by analogs of Stokes data matrices and A. Rudakov's exceptional collections for equivariant quantum differential equations for $\mathbb P^2$.
 
  We also shall  use the $\star-$involution for  the symmetric Laurent polynomial ring 
$$\mathbb C[x^{\pm 1},y^{\pm 1},z^{\pm 1}]^{\mathfrak{S}_3}\simeq \mathbb C[\sigma_1,\sigma_2,\sigma_3^{\pm 1}].$$ Our aim is to simplify the link between generalized  Buchstaber and Kontsevich families of symmetric polynomials. Both families of polynomials considered in our paper can be easily re-written as polynomials in the ring $\mathbb C[\sigma_1,\sigma_2,\sigma_3]$ with generators being elementary symmetric functions in $x,y,z$ (see \eqref{buchstpol} and \eqref{genkontsimm}).

 Consider the symmetric Laurent polynomial ring 
$$\mathbb C[x^{\pm 1},y^{\pm 1},z^{\pm 1}]^{\mathfrak{S}_3}\simeq \mathbb C[\sigma_1,\sigma_2,\sigma_3^{\pm 1}].$$The later isomorphism is given by the evident mapping

$$(x+y+z, xy+yz+zx, xyz)\to (\sigma_1,\sigma_2,\sigma_3).
$$

Both polynomial families have homogeneous degree 2 with respect to weight $(1,1,1)$ and they do not divisible by $\sigma_3$ for generic choice of parameters.
 
 The following lemma is a {\it verbatim} reformulation of the lemma \ref{star} for the case $n=3$:
\begin{lemma}
The involution $\star : \mathbb C[\sigma_1,\sigma_2,\sigma_3^{\pm 1}]\to \mathbb C[\sigma_1,\sigma_2,\sigma_3^{\pm 1}]$ which transform $f\to f^{\star}(x,y,z) = f(\frac{1}{x},\frac{1}{y},\frac{1}{z})$ induces the $\star-$involution on
$\mathbb C[\sigma_1,\sigma_2,\sigma_3^{\pm 1}]:$
$$f^{\star}(\sigma_1,\sigma_2,\sigma_3) = f(\sigma_2/\sigma_3, \sigma_1/\sigma_3,1/\sigma_3)
$$
\end{lemma}
\begin{theorem}\label{th:var}
The following identity holds
$$B_{a_1,a_2,a_3}(\sigma_1,\sigma_2,\sigma_3) = (\sigma_3)^2   D^{\star}_{a_1,a_2,a_3}(\sigma_1,\sigma_2,\sigma_3).
$$ 
\end{theorem}

\section{Conclusion and future perspectives}\label{sec:concl}

In this paper we have started a study of two remarkable polynomial families. These families have some very peculiar common features. It was a challenge to unify them using their common important properties:
\begin{itemize} 
\item both ones belong to the class of so-called discriminantly split symmetric polynomials. Recently, these polynomials were studied in a framework of (generalized) Kovalevsky top integrable systems and their relations to a theory of $2-$valued groups;
\item this common feature describes an interesting class of 2-nd order ODE which includes various equations  such that confluent hypergeometric (2-Bessel et al) and Heun-Lamé equations;
\item products of eigenfunctions for the differential operators entered in the PDEs obey the integral representation with the {\it multiplication kernels}. This classical circle of problems has attracted a lot of attention recently. 
\item one can provide the connection of the multiplication kernels with eigenproblems for generalized shift operators. We show how to get a kernel for generalized shift operators, which started with Kontsevich-Odesskii's description of multiplicative kernels. 
\end{itemize}

We shall address these properties and some important related problems  in the second part of the paper. We are going to stress and investigate  all essential connections between these various features of two polynomial families. The heart of this connection is in the classical product formulas and their modern avatars.

In the second part of the paper, we reconsider the cases when the generalized shift operator with multiplication kernels define a two-valued group structure and when they describe only the structure of a hypergroup. We hope to recover and to extend the Delsarte and Levitan generalized shift operator theories based on the previous results of Buchstaber and Kholodov.
  
We will reconsider the Kovalevskaya-type system integrability for our general polynomial families.
Finally, we shall investigate when the Dragovich classification leads to multiplication kernels which are related to 2-group laws and  when they define only a commutative coalgebra structure on the underlying function spaces.

\end{document}